%% file: Jung.tex
\documentclass[final,hidelinks,onefignum,onetabnum]{template}

\input{shared_info}

\hypersetup{
  colorlinks=false, pdfborder={0 0 1}, linkbordercolor={1 0 0}, citebordercolor={0 1 0} 
}

\begin{document}

\maketitle

\begin{abstract}
    We develop a theory of specular differentiation on real intervals.
    The specular derivative is defined by averaging the angles associated with the forward and backward difference quotients and extends classical differentiation.
    For specularly differentiable and continuous functions, we establish a quasi-mean value theorem and use it to show that continuity of the specular derivative implies $C^1$-regularity.
    Without assuming continuity, we prove that the discontinuities at which the specular derivative is nonzero form an at most countable set; a nowhere-continuous example with identically zero specular derivative shows that this restriction is sharp.
    In the second-order theory of specular differentiation, we impose additional conditions on twice specularly differentiable functions to obtain a class that lies strictly between the classes $C^1$ and $C^2$.
\end{abstract}



\section{Introduction}

Nonsmooth functions can exhibit corners, cusps, and vertical tangents at points where the forward and backward difference quotients fail to converge to a common finite slope.
Generalizations of classical differentiation, such as symmetric differentiation, can assign derivatives at points where the classical derivative does not exist.
Specular differentiation instead uses the geometry of the corresponding one-sided secant lines.
This geometric approach extends classical differentiation while remaining consistent with it.
Moreover, the specular derivative may exist at points where both classical differentiation and other generalized notions of differentiation fail.

\begin{figure}[tbhp]
    \centering
    \subfloat[
        $\partial^+ f(0) = \frac{4}{3}$,
        $\partial^- f(0) = 0$,
        $f^{\sd}(0) = \frac{1}{2}$,
        and $f^{\ast}(0) = \frac{2}{3}$
    ]{%
        \label{fig:specular_finite}%
        \makebox[.43\textwidth][c]{%
            \includegraphics[width=.35\textwidth]
                {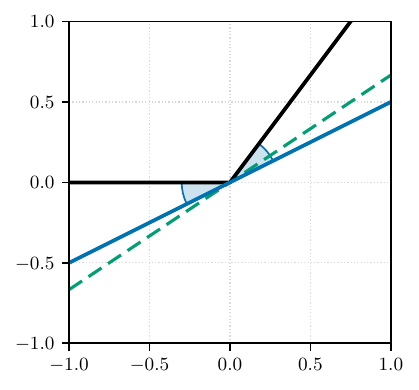}%
        }%
    }
    \hfill
    \subfloat[
        $\partial^+ f(0) = \infty$,
        $\partial^- f(0) = 0$,
        $f^{\sd}(0) = 1$,
        and $f^{\ast}(0)$ does not exist
    ]{%
        \label{fig:specular_infinite}%
        \makebox[.43\textwidth][c]{%
            \includegraphics[width=.35\textwidth]
                {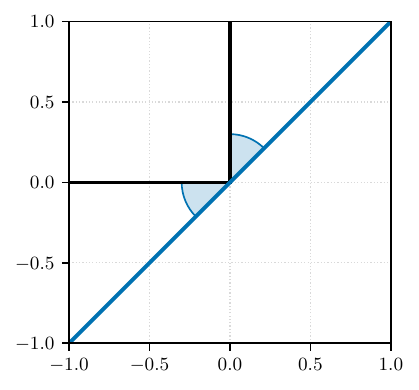}%
        }%
    }
    \caption{Geometric motivation for the specular derivative: black rays (one-sided tangent lines), solid blue line (specular derivative), shaded sectors (equal angles), and green dashed line (symmetric derivative).
    The notation used above is defined in \cref{def:sd,sec:def_and_notation}.}
    \label{fig:motivation}
\end{figure}

The geometric motivation for specular differentiation comes from the law of reflection, according to which the incident and reflected rays make equal angles with the normal to the reflecting surface.
In \cref{fig:motivation}, the black rays represent the two one-sided tangents, while the solid blue line represents their angular bisector and the shaded sectors indicate the equal angles.
In \cref{fig:specular_finite}, the slopes of the solid blue line and the green dashed line represent the specular derivative and the symmetric derivative, respectively.
The latter is obtained by taking the arithmetic mean of the two finite one-sided derivatives.
In \cref{fig:specular_infinite}, one of the one-sided tangents is vertical, whereas the other has a finite slope.
The arithmetic mean is not finite and therefore does not yield a symmetric derivative, whereas the angular bisector has a finite slope and hence yields a specular derivative.

Everywhere specular differentiability implies neither continuity nor Lebesgue measurability, nor does it guarantee the existence of either one-sided derivative, even in the extended-real sense.
This feature is both a limitation and a source of strength: specular differentiability alone does not enforce standard regularity, but it allows the theory to encompass highly irregular functions.
The very breadth of the class of specularly differentiable functions makes it nontrivial to determine how much of classical differential theory remains valid under the present definition.
This question motivates the present paper, which investigates how far classical differential theory, including the Mean Value Theorem and classical regularity results, extends to specular differentiation.

\begin{definition}  \label{def:sd}
    Let $f : I \to \mathbb{R}$ be a function, where $I$ is an open interval in $\mathbb{R}$.
    Let $x \in I$ be a point.
    The \emph{specular derivative} of $f$ at $x$, denoted by $f^{\sd}(x)$, is defined by
    \begin{displaymath}
        f^{\sd}(x) := \lim_{h \searrow 0} \mathcal{A}(f(x + h) - f(x), f(x) - f(x - h), h),
    \end{displaymath}
    where the function $\mathcal{A}:\mathbb{R} \times \mathbb{R} \times (0, \infty) \to \mathbb{R}$ is defined by 
    \begin{displaymath}
        \mathcal{A}(a, b, c) := \frac{a \sqrt{b^2 + c^2} + b \sqrt{a^2 + c^2}}{c \sqrt{a^2 + c^2} + c \sqrt{b^2 + c^2}}.
    \end{displaymath}

    We say $f$ is \emph{specularly differentiable} at $x$ if $f^{\sd}(x)$ exists as a finite real number.
    Moreover, $f$ is said to be \emph{specularly differentiable} on $I$ if $f$ is specularly differentiable at every point in $I$.
\end{definition}

The notion of specular differentiation was first introduced in \cite{2023_Jung}.
The present definition is motivated by \cite[Thm.~2.1.12]{2023_Jung} but differs from the definitions in \cite[Defs.~2.1.5 and~2.2.21]{2023_Jung}.
Throughout this paper, we refer to the derivative defined in \cite{2023_Jung} as the \emph{regular specular derivative}.
This is not merely a terminological distinction: the regular specular derivative is not defined when either one-sided derivative is infinite, whereas the specular derivative introduced here is nevertheless well-defined in some of these cases.
The results of the earlier work cannot simply be transferred to the present setting, and none of them is used as a logical input in the proofs of this paper.
Instead, the corresponding results studied here are formulated using \cref{def:sd} and proved independently within the present framework.

As already observed in \cite{2023_Jung}, classical results on symmetric differentiation motivate several of the questions studied here and provide natural points of comparison.
Specular and symmetric differentiation, however, are incomparable: either derivative may exist at a point where the other does not.
This incomparability stems from two fundamental differences in how the derivatives are constructed.

First, the two notions differ in how they combine the forward and backward difference quotients.
Symmetric differentiation takes their arithmetic mean, whereas specular differentiation takes the arithmetic mean of their arctangents and then converts the resulting angle back into a slope.
In particular, when one difference quotient tends to $+\infty$ or $-\infty$, its arctangent converges to $\frac{\pi}{2}$ or $-\frac{\pi}{2}$, so this limiting angle can still be averaged with the angle associated with the other one-sided quotient; see \cref{fig:specular_infinite}.
More generally, this angular coupling can accommodate certain oscillations for which neither one-sided derivative exists and the symmetric derivative also fails to exist; see \cref{ex:sd_but_not_sy}.
Conversely, arithmetic cancellation may produce a symmetric derivative even when the corresponding angular mean does not converge; see \cref{ex:C0notS1}.

Second, the two notions differ in their dependence on the base-point value $f(x)$.
The value $f(x)$ cancels from the expression defining the symmetric derivative, whereas it appears in both one-sided difference quotients defining the specular derivative.
Thus, changing $f(x)$ alone does not affect the existence or value of the symmetric derivative at $x$.
By contrast, the same change may alter the value of the specular derivative or prevent it from existing.

For background and broader accounts of symmetric differentiation, see Larson \cite{1981_Larson,19831984_Larson} and the comprehensive monograph by Thomson \cite{1994_Thomson_BOOK}.
A particularly relevant point of comparison is the Quasi-Mean Value Theorem of Aull \cite{1967_Aull} for functions that are continuous on a closed interval and symmetrically differentiable on its interior.
For further discussion and related versions of this theorem, see \cite{19911992_Evans,1998_Sahoo_BOOK,2011_Sahoo}.

Symmetric and classical differentiation exhibit markedly different regularity behavior.
This contrast is illustrated by the examples in \cite{1977_Foran,1985_Ponomarev,19871988_Uher,19881989_Freiling}.
Several strong regularity results are known for symmetrically differentiable functions; see \cite{1978_Belna,19821983_Uher,1983_Larson}.
Another aspect of the theory concerns the structure of discontinuity sets; see \cite{1933_Charzynski,1990_Freiling}.
We later compare these results with the corresponding conclusions for specular differentiation.
These comparisons require separate arguments and may involve different hypotheses.

\subsection{Definitions and notation}
\label{sec:def_and_notation}

Let $\overline{\mathbb{R}} := \mathbb{R} \cup \{-\infty, +\infty\}$ denote the extended real line.
Let $\mathcal{L}^1$ denote Lebesgue measure on $\mathbb{R}$.
For a set $E \subset \mathbb{R}$, let $\chi_E$ denote the characteristic function of $E$.
Let $\sgn:\mathbb{R}\to\mathbb{R}$ denote the sign function, with $\sgn(0):=0$.

Let $I$ be an open interval in $\mathbb{R}$, and $E$ be a subset of $I$.
The set $E$ is called a \emph{$G_\delta$ set} in $I$ if there exist open subsets $U_n$ of $I$ such that $E=\bigcap_{n=1}^{\infty}U_n$.
The set $E$ is called \emph{nowhere dense} in $I$ if the interior in $I$ of its closure in $I$ is empty.
The set $E$ is called \emph{meager} in $I$ if it is a countable union of nowhere dense subsets of $I$.
The set $E$ is called \emph{residual} in $I$ if $I\setminus E$ is meager in $I$.

Let $C^0(I)$ denote the space of continuous real-valued functions on $I$.
We denote by $C^1(I)$ the space of continuously differentiable real-valued functions on $I$.
Let $B^1(I)$ denote the class of Baire class $1$ functions.

We introduce the following classes of specularly differentiable functions.

\begin{definition}
    We define
    \begin{align*}
        S^0 (I)
        &:= \left\{ f : I \to \mathbb{R} \, \middle| \, f \text{ is specularly differentiable on } I \right\},   \\
        S^1 (I)
        &:= C^0(I) \cap S^0(I).
    \end{align*}
\end{definition}

\noindent Let $f : I \to \mathbb{R}$ be a function.
We denote the sets of continuity and discontinuity points of $f$ by
\begin{align*}
    \Cont(f)
    &:= \left\{ x \in I \, \middle| \, f \text{ is continuous at } x \right\} ,  \\
    \Disc(f)
    &:= \left\{ x \in I \, \middle| \, f \text{ is discontinuous at } x \right\}.
\end{align*}

\noindent Let $x \in I$ be a point.
For $h>0$ such that $x-h, x+h\in I$, we denote the \emph{right} and \emph{left quotients} of $f$ at $x$ by
\begin{displaymath}
    q^+_h f(x) := \frac{f(x + h) - f(x)}{h} 
    \qquad\text{and}\qquad
    q^-_h f(x) := \frac{f(x) - f(x - h)}{h},
\end{displaymath}
respectively.
The \emph{right} and \emph{left derivatives} of $f$ at $x$ are denoted by 
\begin{displaymath}
    \partial^+ f(x) :=  \lim_{h \searrow 0} q^+_h f(x) 
    \qquad\text{and}\qquad
    \partial^- f(x) :=  \lim_{h \searrow 0} q^-_h f(x),
\end{displaymath}
respectively.
These derivatives may take values in $\overline{\mathbb{R}}$.
Collectively, $\partial^+f(x)$ and $\partial^-f(x)$ are called the \emph{one-sided derivatives} of $f$ at $x$.
The \emph{symmetric derivative} of $f$ at $x$ is defined by 
\begin{displaymath}
    f^{\ast}(x) 
    := 
    \lim_{h \to 0} \frac{f(x + h) - f(x - h)}{2h} 
    = 
    \lim_{h \searrow 0} \frac{q^+_h f(x) + q^-_h f(x)}{2}
\end{displaymath}
if the limit exists as a finite real number.

We introduce auxiliary functions which are particularly useful in formulating specular derivatives and their estimates.
The parameters $\alpha$ and $\beta$ below are intended to represent right and left derivatives.

\begin{definition}  \label{def:BC}
    First, define the function $\mathcal{B} : \overline{\mathbb{R}}^2 \to \overline{\mathbb{R}}$ by 
    \begin{equation}    \label{def:B}
        \mathcal{B}(\alpha, \beta) :=  \displaystyle \tan \left( \frac{1}{2} \arctan ( \alpha ) + \frac{1}{2} \arctan (\beta) \right) ,
    \end{equation}
    where we adopt the conventions $\arctan ( \pm \infty):= \pm \frac{\pi}{2}$ and $\tan \left( \pm \frac{\pi}{2}\right):= \pm \infty$.
    Second, define the function $\mathcal{C}: \overline{\mathbb{R}}^2 \to \overline{\mathbb{R}}$ by
    \begin{displaymath} 
        \mathcal{C}(\alpha, \beta) := 
        \begin{cases}
            \left( \frac{\alpha}{\sqrt{1 + \alpha^2}} + \frac{\beta}{\sqrt{1 + \beta^2}} \right) \left( \frac{1}{\sqrt{1 + \alpha^2}} + \frac{1}{\sqrt{1 + \beta^2}} \right)^{-1}
            & \text{if $\alpha, \beta \in \mathbb{R}$,}  \\
            \alpha \pm \sqrt{1 + \alpha^2} 
            & \text{if $\alpha \in \mathbb{R}$, $\beta = \pm \infty$,}   \\
            \beta \pm \sqrt{1 + \beta^2} 
            & \text{if $\alpha = \pm \infty$, $\beta \in \mathbb{R}$,}  \\
            0 
            & \text{if $(\alpha,\beta)=(\pm\infty,\mp\infty)$,} \\
            \pm \infty 
            & \text{if $(\alpha,\beta)=(\pm\infty,\pm\infty)$}.
        \end{cases}
    \end{displaymath}
\end{definition}

\subsection{Main results}

We first analyze the relationship between specular and symmetric differentiation, identifying both common features and fundamental differences.
In the presence of a one-sided derivative, however, the specular derivative admits a representation in terms of the two one-sided derivatives.
More precisely, if $f$ is specularly differentiable at a point and at least one of the one-sided derivatives $\partial^+f$ and $\partial^-f$ exists there in $\overline{\mathbb{R}}$, then both exist there and their angular mean equals the specular derivative at that point; see \cref{thm:repr_sd}.

Let $I$ be an open interval.
If both one-sided derivatives of a function $f:I\to\mathbb{R}$ exist as finite real numbers at every point of the open interval $I$, then $f\in S^1(I)$; see \cref{prop:finite_one_sided}.
Consequently, every real-valued convex function $f$ on $I$ belongs to $S^1(I)$, and $f^{\sd}(x)\in\partial f(x)$ for every $x\in I$, where $\partial f(x)$ denotes the subdifferential of $f$ at $x$ in the sense of convex analysis; see \cref{cor:convex_sd_subgradient}.

For continuous specularly differentiable functions, our central result is the Quasi-Mean Value Theorem; see \cref{thm:q-MVT}.
This theorem bounds every secant slope from below and above by values of the specular derivative attained in the interior of the interval.
It retains the essential order information of the classical Mean Value Theorem without requiring the secant slope to equal the specular derivative at a single point.
As a consequence, the sign of the specular derivative characterizes monotonicity, and the specular derivative vanishes identically if and only if the function is constant; see \cref{cor:monotonicity_and_sd}.

We next establish a countability result without assuming continuity or measurability.
If $f\in S^0(I)$, then
\begin{displaymath}
    \left\{
        x\in\Disc(f)
        \,\middle|\,
        f^{\sd}(x)\neq0
    \right\}
\end{displaymath}
is at most countable; see
\cref{thm:countable_nonzero_sd_discontinuities}.
The restriction to points at which the specular derivative is nonzero is essential.
Indeed, \cref{ex:nowhere_continuous} gives a nowhere-continuous function $f:\mathbb{R}\to\mathbb{R}$ that is not Lebesgue measurable and satisfies $f^{\sd} \equiv 0$.

We next establish regularity properties of specularly differentiable functions and their specular derivatives.
Let $f\in S^1(I)$.
Then the specular derivative $f^{\sd}$ is of Baire class $1$; see \cref{prop:Baire_one}.
Consequently, $\Cont(f^{\sd})$ is a dense $G_\delta$ subset of $I$.
Together with the Quasi-Mean Value Theorem, this implies that $f$ is classically differentiable at every point of $\Cont(f^{\sd})$ and that $f'=f^{\sd}$ on $\Cont(f^{\sd})$; see \cref{thm:residual_differentiability}.
Moreover, if $f^{\sd}$ is bounded on $I$, then $f$ is Lipschitz continuous and $f'=f^{\sd}$ $\mathcal{L}^1$-a.e. on $I$; see \cref{prop:cont_bounded_sd_impy_Lip_cont}.
If $f^{\sd}$ is continuous on $I$, then $f\in C^1(I)$ and $f'=f^{\sd}$ on $I$; see \cref{thm:cont_sd_implies_d}.

Finally, we identify additional continuity conditions for twice specularly differentiable functions and use them to define the class $S^2(I)$; see \cref{def:S2}.
We then prove the strict inclusions
\begin{displaymath}
    C^2(I)
    \subsetneq
    S^2(I)
    \subsetneq
    C^1(I)
    \subsetneq
    S^1(I)
    \subsetneq
    C^0(I);
\end{displaymath}
see \cref{thm:chain_regularity}.
The requirement that the specular derivative be continuous is essential: we construct a continuous, twice specularly differentiable function that is not classically differentiable; see \cref{ex:notS2notC1}.

\subsection{Organization}

The organization of this paper is as follows.
In \cref{sec:analysis}, we establish equivalent representations of the specular derivative and compare it with the symmetric derivative.
In \cref{sec:Q-MVT}, we prove Quasi-Rolle's Theorem and the Quasi-Mean Value Theorem, and derive characterizations of monotonicity and constancy in terms of the sign of the specular derivative.
In \cref{sec:discontinuities}, we investigate the discontinuities of specularly differentiable functions.
In \cref{sec:regularity}, we investigate regularity properties of both specularly differentiable functions and their specular derivatives.
In \cref{sec:2nd_order}, we introduce second-order specular differentiability, compare it with classical smoothness, and show that twice specular differentiability alone does not imply classical differentiability.
Finally, \cref{apx:A_and_B} analyzes the auxiliary functions $\mathcal{A}$, $\mathcal{B}$, and $\mathcal{C}$; \cref{apx:Leicester} states and proves Young's Leicester theorem.

\section{Analysis of specular differentiation}   
\label{sec:analysis}

The specular derivative can be represented in two equivalent forms by means of the following relationship among the three auxiliary functions.

\begin{lemma} \label{lem:ABC}
    For every $(a, b, c) \in \mathbb{R} \times \mathbb{R} \times (0, \infty)$, we have
    \begin{displaymath}
        \mathcal{A}(a, b, c)
        = \mathcal{B}\left( \frac{a}{c}, \frac{b}{c} \right) 
        = \mathcal{C} \left( \frac{a}{c}, \frac{b}{c} \right).
    \end{displaymath}
    Furthermore, for every $(\alpha,\beta)\in\overline{\mathbb{R}}^2$, we have $\mathcal{B}(\alpha,\beta)=\mathcal{C}(\alpha,\beta)$.
\end{lemma}

\begin{proof}
    See \cref{apx:A_and_B}.
\end{proof}

Let $f: I \to \mathbb{R}$ be a function, where $I$ is an open interval in $\mathbb{R}$.
Let $x \in I$ be a point. 
By \cref{lem:ABC}, the defining expression for the specular derivative can be written in the following two equivalent forms:
\begin{align}
    f^{\sd}(x) 
    &=\lim_{h \searrow 0} \mathcal{B} \left( q^+_h f(x), q^-_h f(x) \right) \label{eq:1st_form_sd} \\
    &= \lim_{h \searrow 0} \mathcal{C} \left( q^+_h f(x), q^-_h f(x) \right). \label{eq:2nd_form_sd}
\end{align}

\begin{remark}
    Specular differentiation generalizes classical differentiation.
    If $f$ is classically differentiable at $x$, then both difference quotients in \eqref{eq:1st_form_sd} converge to $f'(x)$.
    Hence, by \cref{lem:ABC} and \cref{lem:A}~\ref{lem:A-1},
    \begin{displaymath}
        f'(x)
        =
        \mathcal{B}\bigl(f'(x),f'(x)\bigr)
        =
        f^{\sd}(x).
    \end{displaymath}
\end{remark}

\subsection{Specular and symmetric derivatives}

Specular and symmetric derivatives share some properties but differ in others.
The first common feature is as follows.

\begin{example} \label{ex:even}
    The specular derivative of any locally symmetric function is zero.
    More precisely, if $f(x - h) = f(x + h)$ for sufficiently small $h > 0$, then $f^{\sd}(x) = 0$.
    This directly follows from the formula \cref{eq:1st_form_sd}.
    This property is shared with symmetric derivatives: the symmetric derivative of any locally symmetric function is also zero.

    A trivial example is the characteristic function $\chi_{\left\{ 0 \right\}} : \mathbb{R} \to \mathbb{R}$ defined by 
    \begin{displaymath}
        \chi_{\left\{ 0 \right\}} (x) =
        \begin{cases}
            0    &    \text{if } x \neq 0,    \\
            1    &    \text{if } x = 0.
        \end{cases}
    \end{displaymath}
    This function is discontinuous at $x=0$.
    Then, $(\chi_{\left\{ 0 \right\}})^{\sd}(0) = 0 = (\chi_{\left\{ 0 \right\}})^{\ast}(0)$.

    A nontrivial example is the even function $f:\mathbb{R} \to \mathbb{R}$ defined by 
    \begin{equation} \label{eq:osc_zero_sd}
        f(x) =
        \begin{cases}
            \displaystyle  x \sin \frac{1}{x}    &    \text{if } x \neq 0,    \\[2mm]
            0    &    \text{if } x = 0.
        \end{cases}
    \end{equation} 
    Then, we have $f^{\sd}(0) = 0 = f^{\ast}(0)$.
    Note that $\partial^+ f(0)$ and $\partial^- f(0)$ do not exist.
\end{example}

Second, the Quasi-Mean Value Theorem for specular derivatives (\cref{thm:q-MVT}) holds, as it does for symmetric derivatives (\cite[Thm.~1]{1967_Aull}).
Some regularity byproducts of the Quasi-Mean Value Theorem in the specular sense parallel those for symmetric derivatives; we defer this discussion to \cref{sec:Q-MVT,sec:regularity}, as it requires further analysis of specular derivatives.

The third common feature of specular and symmetric derivatives is that at a local extremum within an open interval, neither derivative need be zero.

\begin{example}
    The function $f(x) = x - |2x|$, for $x \in \mathbb{R}$, has a maximum at $x=0$.
    Then, $f^{\sd}(0) = \sqrt{5} - 2$ and $f^{\ast}(0) = 1$.
    This example is adapted from \cite[Sect.~4]{1967_Aull}.
\end{example}

However, a bound on the specular derivative can still be obtained.

\begin{remark}
    Let $f : I \to \mathbb{R}$ be specularly differentiable on $I$, where $I$ is an open interval in $\mathbb{R}$.
    If $x^{\ast} \in I$ is a local extremum of $f$, then $\left\vert f^{\sd}(x^{\ast}) \right\vert \leq 1$.
    Indeed, suppose that $x^\ast$ is a local minimum point of $f$.
    For every sufficiently small $h>0$, we have
    \begin{displaymath}
        0 \leq \frac{1}{2}\arctan \left( \frac{f(x^{\ast} + h) - f(x^{\ast})}{h} \right) \leq \frac{\pi}{4}
    \end{displaymath}   
    and 
    \begin{displaymath}
        -\frac{\pi}{4} \leq \frac{1}{2} \arctan \left( \frac{f(x^{\ast}) - f(x^{\ast} - h)}{h} \right) \leq 0.
    \end{displaymath}
    Adding these inequalities, applying the increasing tangent function, and taking the limit as $h \searrow 0$ in \eqref{eq:1st_form_sd}, we obtain the desired inequality.
    The conclusion for a local maximum point follows similarly.
\end{remark}

There are, of course, differences between specular and symmetric derivatives.

\begin{remark} \label{rmk:dependence_on_base_value}
    A fundamental distinction between the symmetric and specular derivatives is the role of the base-point value $f(x)$.
    The expression defining the symmetric derivative is the arithmetic mean of the forward and backward difference quotients:
    \begin{displaymath}
        \frac{f(x+h)-f(x-h)}{2h}
        =
        \frac{q_h^+f(x)+q_h^-f(x)}{2}.
    \end{displaymath}
    Although each difference quotient separately involves $f(x)$, the two occurrences of $f(x)$ cancel in their arithmetic mean.
    Hence, the existence and value of the symmetric derivative at $x$ depend only on the restriction of $f$ to a punctured neighborhood of $x$ and are unchanged if only $f(x)$ is modified.

    By contrast, the corresponding expression for the specular derivative is
    \begin{displaymath}
        \mathcal{B} \left( \frac{f(x + h) - f(x)}{h}, \frac{f(x) - f(x - h)}{h} \right)
        = \mathcal{B}\bigl(q_h^+f(x),q_h^-f(x)\bigr),
    \end{displaymath}
    in which the dependence on $f(x)$ does not generally cancel.
    Thus, the specular derivative is sensitive to the value of the function at the base point: changing $f(x)$ while leaving the function unchanged on a punctured neighborhood of $x$ can alter the value or even destroy the existence of the specular derivative.
    This phenomenon is illustrated in \cref{ex:l.s.c._func,ex:disc_func}.
\end{remark}

Specular and symmetric differentiability do not imply one another, as the following examples show.

\begin{example} \label{ex:sd_but_not_sy}
    Specular differentiability does not imply symmetric differentiability.
    Consider the non-even function $f:\mathbb{R} \to \mathbb{R}$ defined by
    \begin{displaymath}
        f(x) =
        \begin{cases}
            \displaystyle x \left( 2 + \sin \frac{1}{x} \right)    &    \text{if } x > 0 ,    \\[2mm]
            0    &    \text{if } x = 0,    \\
            \displaystyle x \left( 2 - \sin \frac{1}{x} \right)^{-1}    &    \text{if } x < 0.
        \end{cases}
    \end{displaymath}
    The identity $\arctan(t) + \arctan(t^{-1}) = \frac{\pi}{2}$ for all $t > 0$ gives
    \begin{displaymath}
        \arctan\left( \frac{f(0 + h) - f(0)}{h} \right) + \arctan\left( \frac{f(0) - f(0 - h)}{h} \right) = \frac{\pi}{2}.
    \end{displaymath}
    Therefore, the formula \cref{eq:1st_form_sd} yields $f^{\sd}(0) = 1$, but $\partial^+ f(0)$ and $\partial^- f(0)$ do not exist.
    However, the symmetric derivative $f^{\ast}(0)$ does not exist, since 
    \begin{displaymath}
        \frac{f(h) - f(-h)}{2h} = 1 + \frac{1}{2} \sin \frac{1}{h} + \left( 4 + 2 \sin \frac{1}{h} \right)^{-1},
    \end{displaymath}
    which does not converge as $h \searrow 0$.
\end{example}

\begin{example}[A function in $C^0$ but not in $S^1$] \label{ex:C0notS1}
    Symmetric differentiability does not imply specular differentiability.
    Define the function $f : \mathbb{R} \to \mathbb{R}$ by
    \begin{displaymath}
        f(x):=
        \begin{cases}
            \displaystyle x\sin \frac{1}{x}, & x > 0,\\
            0, & x=0,\\
            \displaystyle x\left(1+\sin \frac{1}{x}\right), & x < 0.
        \end{cases}
    \end{displaymath}
    For $h>0$, we have $q^+_h f(0) = \sin \frac{1}{h}$ and $q^-_h f(0) = 1 - \sin \frac{1}{h}$.
    Neither $\partial^+f(0)$ nor $\partial^-f(0)$ exists.
    On the one hand, $f$ is symmetrically differentiable at $x=0$, with $f^{\ast}(0)=\frac{1}{2}$.
    On the other hand, $f$ is not specularly differentiable at $x=0$.
    Indeed, for each $n\in\mathbb{N}$, let
    \begin{displaymath}
        h_n:=\frac{1}{n\pi}
        \qquad\text{and}\qquad
        k_n:=\frac{6}{\pi+12n\pi}.
    \end{displaymath}
    Then both sequences converge to zero.
    For every $n\geq1$,
    \begin{displaymath}
        \mathcal{B}\left(q_{h_n}^+f(0),q_{h_n}^-f(0)\right)
        =
        \mathcal{B}(0,1)
        =
        \sqrt{2}-1,
    \end{displaymath}
    whereas
    \begin{displaymath}
        \mathcal{B}\left(q_{k_n}^+f(0),q_{k_n}^-f(0)\right)
        =
        \mathcal{B}\left(\frac{1}{2},\frac{1}{2}\right)
        =
        \frac{1}{2}.
    \end{displaymath}
    Thus, the expression in \eqref{eq:1st_form_sd} has two distinct subsequential limits, and hence $f^{\sd}(0)$ does not exist.
    Since $f$ is continuous on $\mathbb{R}$, this shows that $f\in C^0(\mathbb{R})$ but $f\notin S^1(\mathbb{R})$.
\end{example}

A further distinction appears when the one-sided derivatives exist in $\overline{\mathbb{R}}$.
If one is finite and the other is infinite, then the specular derivative exists as a finite real number, whereas the symmetric derivative does not.
The following example illustrates this case.

\begin{example}[A function in $S^1$ but not in $C^1$] \label{ex:S1notC1}
    Consider the function $f : \mathbb{R} \to \mathbb{R}$ defined by 
    \begin{displaymath}
        f(x)
        := 
        \begin{cases}
            \sqrt{x}        &    \mbox{if } x \geq 0,    \\
            -x        &    \mbox{if } x < 0.    
        \end{cases}
    \end{displaymath}
    A direct calculation shows that $f^{\sd}(0) = \sqrt{2} -1$. 
    Since $f$ is continuous on $\mathbb{R}$ and classically differentiable away from $0$, it follows that $f \in S^1(\mathbb{R})$.
    However, $f'(0)$ does not exist, and hence $f\notin C^1(\mathbb{R})$.
    In fact, even the symmetric derivative $f^{\ast}(0)$ does not exist.
\end{example}

Symmetric differentiation is linear.
In contrast, specular differentiation need not be homogeneous; moreover, the class $S^0$ need not be closed under addition.

\begin{example}[Failure of linearity and closure under addition] \label{ex:sd_not_preserved_addition}
    Consider the functions $f, g:\mathbb{R} \to \mathbb{R}$ defined by 
    \begin{displaymath}
        f(x) =
        \begin{cases}
            0    &    \text{if } x \in (-\infty, 0],    \\
            1    &    \text{if } x \in (0, \infty)      \\
        \end{cases}
        \qquad\text{and}\qquad
        g(x) = 
        \begin{cases}
            0    &    \text{if } x \in (-\infty, 0),    \\
            1    &    \text{if } x \in [0, \infty) .
            \end{cases}
    \end{displaymath}
    Then, $f$ is lower semicontinuous and $g$ is upper semicontinuous on $\mathbb{R}$.
    Both functions are locally constant away from zero, and a direct calculation gives $f^{\sd}(0) = 1 = g^{\sd}(0)$.
    Moreover,
    \begin{displaymath}
        (2f)^{\sd}(0)=1\neq2=2f^{\sd}(0),
    \end{displaymath}
    which shows that specular differentiation is not homogeneous and therefore not linear.

    On the other hand, both difference quotients of $(f + g)$ at $x = 0$ are equal to $\frac{1}{h}$ for $h>0$.
    Therefore, $(f + g)^{\sd}(0)$ does not exist as a finite real number.
    Thus, $f, g\in S^0(\mathbb{R})$, but $f + g \notin S^0(\mathbb{R})$.
\end{example}

For comparison, if $f$ is continuous at $x$ and symmetrically differentiable at $x$, and $\varphi$ is continuously differentiable in a neighborhood of $f(x)$, then
\begin{displaymath}
    (\varphi\circ f)^{\ast}(x)
    =
    \varphi'(f(x))f^{\ast}(x).
\end{displaymath}
The corresponding chain rule may fail for specular differentiation, even for a continuous inner function and a smooth outer function.

\begin{example}[Failure of the chain rule] \label{ex:failure_chain_rule}
    Let $f$ be the function in \cref{ex:S1notC1}, and define
    $\varphi:\mathbb{R}\to\mathbb{R}$ by $\varphi(t):=t^2$.
    For $h>0$, we have
    \begin{displaymath}
        \frac{(\varphi\circ f)(h)-(\varphi\circ f)(0)}{h}=1
        \qquad\text{and}\qquad
        \frac{(\varphi\circ f)(0)-(\varphi\circ f)(-h)}{h}=-h.
    \end{displaymath}
    Hence,
    \begin{displaymath}
        (\varphi\circ f)^{\sd}(0)
        =
        \mathcal{B}(1,0)
        =
        \sqrt{2}-1
        \neq 
        0
        =
        \varphi'(f(0))f^{\sd}(0).
    \end{displaymath}
    Thus, the usual chain rule fails, even though $f\in S^1(\mathbb{R})$ and $\varphi\in C^\infty(\mathbb{R})$.
\end{example}

Specular and symmetric differentiability also exhibit different regularity properties, both for the underlying functions and for their corresponding derivatives.
We investigate these differences in \cref{sec:regularity}.

\subsection{Representations of specular derivatives in terms of one-sided derivatives}   

The existence of the specular derivative and one of the one-sided derivatives (possibly infinite) implies the existence of the other one-sided derivative.

\begin{lemma}   \label{lem:existence_one_sided_d}
    Let $f: I \to \mathbb{R}$ be a function, where $I$ is an open interval in $\mathbb{R}$.
    Fix a point $x \in I$.
    Suppose that the specular derivative $f^{\sd}(x)$ exists as a real number.
    Then, the following statements hold.
    \begin{enumerate}[label=\upshape(\alph*)]
        \item \label{lem:existence_one_sided_d-a} If $\partial^+ f(x)$ exists as an extended real number, then $\partial^- f(x)$ exists as an extended real number. 
        \item \label{lem:existence_one_sided_d-b} If $\partial^- f(x)$ exists as an extended real number, then $\partial^+ f(x)$ exists as an extended real number. 
        \item \label{lem:existence_one_sided_d-c} If both one-sided derivatives exist in $\overline{\mathbb{R}}$, then they cannot both equal $+\infty$ or both equal $-\infty$. 
    \end{enumerate}
\end{lemma}

\begin{proof}
    For all sufficiently small $h>0$, set $\alpha_h := q_h^+f(x)$ and $\beta_h := q_h^-f(x)$.
    Also, write $\gamma:=f^{\sd}(x)$, and, whenever they exist in $\overline{\mathbb{R}}$, denote the one-sided derivatives by $\alpha:=\partial^+f(x)$ and $\beta:=\partial^-f(x)$.

    First, we prove \cref{lem:existence_one_sided_d-a}.
    Let $\varepsilon >0$ be arbitrary.
    By the existence of $\gamma$ and the formula \eqref{eq:1st_form_sd}, we find that 
    \begin{displaymath}
        2 \arctan \gamma = \lim_{h \searrow 0} \left( \arctan \alpha_h + \arctan \beta_h \right).
    \end{displaymath}
    Thus, there exists $\delta_1(\varepsilon) > 0$ such that, if $h \in (0, \delta_1)$, then
    \begin{equation}    \label{ineq:repr_of_sd-1}
        \left\vert \, \arctan \alpha_h + \arctan \beta_h - 2 \arctan \gamma \, \right\vert < \frac{\varepsilon}{2}.
    \end{equation}
    From the existence of $\alpha$ with the convention $\arctan(\pm \infty) = \pm \frac{\pi}{2}$, there exists $\delta_2(\varepsilon) > 0$ such that, if $h \in (0, \delta_2)$, then 
    \begin{equation}   \label{ineq:repr_of_sd-2}
        \left\vert \, \arctan \alpha_h - \arctan \alpha \right\vert < \frac{\varepsilon}{2}.
    \end{equation}
    If $0 < h < \min\left\{ \delta_1, \delta_2 \right\}$, then inequalities \cref{ineq:repr_of_sd-1} and \cref{ineq:repr_of_sd-2} imply that
    \begin{align*}
        & \left\vert \, \arctan \beta_h - (2 \arctan \gamma - \arctan \alpha ) \right\vert \\
        \leq & \left\vert \, \arctan \alpha_h + \arctan \beta_h - 2 \arctan \gamma \right\vert + \left\vert \, \arctan \alpha_h - \arctan \alpha \right\vert \\
        <& \, \varepsilon.
    \end{align*}
    Since $\varepsilon>0$ is arbitrary, we have
    \begin{displaymath}
        \lim_{h\searrow0} \arctan \beta_h
        =
        2\arctan \gamma - \arctan \alpha .
    \end{displaymath}
    The limiting angle belongs to $\big[-\frac{\pi}{2},\frac{\pi}{2}\big]$.
    By the conventions in \eqref{def:B}, the extended arctangent $\arctan : \overline{\mathbb{R}} \to \big[ -\frac{\pi}{2},\frac{\pi}{2} \big]$ is a homeomorphism whose inverse is the extended tangent function.
    Consequently,
    \begin{displaymath}
        \lim_{h \searrow 0} \beta_h
        =
        \tan\left( 2\arctan \gamma - \arctan \alpha \right) 
        \in \overline{\mathbb{R}}.
    \end{displaymath}
    Thus, $\partial^-f(x)$ exists as an extended real number.
    \Cref{lem:existence_one_sided_d-b} can be proved similarly.

    Next, to show \cref{lem:existence_one_sided_d-c}, suppose to the contrary that $\alpha = \infty = \beta$.
    By the existence of $\gamma$ and the formula \eqref{eq:1st_form_sd}, there exists $\delta_3 > 0$ such that, if $h \in (0, \delta_3)$, then $\left\vert \, \mathcal{B} \left( \alpha_h, \beta_h \right) - \gamma \, \right\vert  < 1$, which implies that  
    \begin{equation} \label{ineq:repr_of_sd-3}
        \arctan \alpha_h + \arctan \beta_h < 2\arctan (1 + \gamma).
    \end{equation}
    Since $\alpha = \infty = \beta$, there exists $\delta_4 > 0$ such that, if $h \in (0, \delta_4)$, then $\alpha_h > 2 + |\gamma|$ and $\beta_h > 2 + |\gamma|$.
    Combining these inequalities gives
    \begin{displaymath}
        \arctan \alpha_h + \arctan \beta_h
        > 2\arctan (2 + |\gamma|)
        > 2\arctan (1 + \gamma),
    \end{displaymath}
    which contradicts the inequality \eqref{ineq:repr_of_sd-3} when $0 < h < \min\{\delta_3, \delta_4\}$.
    Therefore, the case $\alpha = \infty = \beta$ is impossible.
    The other case $\alpha = -\infty = \beta$ can be proved similarly.
\end{proof}

\begin{remark}
    The additional assumption that at least one of the one-sided derivatives exists cannot be omitted from \cref{lem:existence_one_sided_d}.
    In other words, the existence of a specular derivative does not imply the existence of any one-sided derivative.
    The even function defined in \cref{eq:osc_zero_sd} and the non-even function defined in \cref{ex:sd_but_not_sy} serve as counterexamples.
\end{remark}

If a function is specularly differentiable at a point and has either a left or a right derivative there, possibly infinite, then the other one-sided derivative also exists, and the specular derivative is determined by the two one-sided derivatives.

\begin{theorem} \label{thm:repr_sd}
    Let $f: I \to \mathbb{R}$ be a function, where $I$ is an open interval in $\mathbb{R}$.
    Fix a point $x \in I$.
    Assume that $f^{\sd}(x)$ exists and that at least one of the one-sided derivatives $\partial^+f(x)$ and $\partial^-f(x)$
    exists as an extended real number.
    Then both one-sided derivatives exist in $\overline{\mathbb{R}}$, and
    \begin{align}
        f^{\sd}(x) 
        &= \mathcal{B}\left( \partial^+ f(x), \partial^- f(x) \right), \label{eq:1st_form_sd_wo_limit}   \\
        &= \mathcal{C}\left( \partial^+ f(x), \partial^- f(x) \right). \label{eq:2nd_form_sd_wo_limit}
    \end{align}
    In particular, the common value is a real number.
\end{theorem}

\begin{proof}
    By \cref{lem:existence_one_sided_d-a,lem:existence_one_sided_d-b} of \cref{lem:existence_one_sided_d}, both one-sided derivatives exist in $\overline{\mathbb{R}}$.
    Moreover, \cref{lem:existence_one_sided_d-c} of \cref{lem:existence_one_sided_d} implies that
    \begin{displaymath}
        \frac{1}{2} \left( \arctan(\partial^+f(x)) + \arctan(\partial^-f(x)) \right) \in \left(-\frac{\pi}{2},\frac{\pi}{2}\right).
    \end{displaymath}
    Hence, taking $h \searrow 0$ in \eqref{eq:1st_form_sd} yields the first equality \eqref{eq:1st_form_sd_wo_limit}.
    The second equality \eqref{eq:2nd_form_sd_wo_limit} follows from \cref{lem:ABC}.
    The common value is finite since the averaged angle belongs to $\left(-\frac{\pi}{2},\frac{\pi}{2}\right)$.
\end{proof}

When both one-sided derivatives are finite, the formula \eqref{eq:2nd_form_sd_wo_limit} agrees with the corresponding result for regular specular derivatives in \cite[Prop.~2.1.15]{2023_Jung}.
When exactly one of the one-sided derivatives is infinite, the formula agrees with the limiting identities in \cref{lem:A}~\ref{lem:A-13}.

Once bounds on the one-sided derivatives are known, the specular derivative can be estimated accordingly.

\begin{corollary} \label{cor:estimate_of_sd}
    Let $f: I \to \mathbb{R}$ be a function, where $I$ is an open interval in $\mathbb{R}$.
    Fix a point $x \in I$.
    Assume that the specular derivative $f^{\sd}(x)$ exists and that at least one of the one-sided derivatives $\partial^{\pm}f(x)$ exists as an extended real number.
    Then, the following statements hold.
    \begin{enumerate}[label=\upshape(\alph*)]
        \item \label{cor:estimate_of_sd-1} Either $\partial^- f(x) \leq f^{\sd}(x) \leq \partial^+ f(x)$ or $\partial^+ f(x) \leq  f^{\sd}(x) \leq \partial^- f(x)$.
        \item \label{cor:estimate_of_sd-2} If $m_1, m_2 \in \overline{\mathbb{R}}$ with $m_1 \leq \partial^+ f(x)$ and $m_2 \leq \partial^- f(x)$, then 
        \begin{displaymath}
            \min\left\{ m_1, m_2 \right\} \leq f^{\sd}(x).
        \end{displaymath}
        \item \label{cor:estimate_of_sd-3} If $M_1, M_2 \in \overline{\mathbb{R}}$ with $\partial^+ f(x) \leq M_1$ and $\partial^- f(x) \leq M_2$, then 
        \begin{displaymath}
            f^{\sd}(x) \leq \max\left\{ M_1, M_2 \right\}.
        \end{displaymath}
    \end{enumerate}
\end{corollary}

\begin{proof}
    The existence of $\partial^+ f(x)$ and $\partial^- f(x)$ follows from \cref{lem:existence_one_sided_d}.
    Define $m := \min \left\{ \partial^+ f(x), \partial^- f(x) \right\}$ and $M := \max \left\{ \partial^+ f(x), \partial^- f(x) \right\}$.
    By \cref{thm:repr_sd}, we have $f^{\sd}(x) = \mathcal{B} \left( \partial^+ f(x), \partial^- f(x) \right)$.
    Thus, by \cref{lem:monotonicity_of_B} and the identity $\mathcal{B}(t,t)=t$ for every $t\in\overline{\mathbb{R}}$, which follows directly from \eqref{def:B}, we obtain
    \begin{equation} \label{ineq:estimate_of_sd-1}
        m 
        = 
        \mathcal{B} (m, m) 
        \leq 
        \mathcal{B} \left( \partial^+ f(x), \partial^- f(x) \right) 
        \leq 
        \mathcal{B} (M, M) 
        = 
        M,
    \end{equation}
    which proves \cref{cor:estimate_of_sd-1}.
    To prove \cref{cor:estimate_of_sd-2}, let $m_1, m_2 \in \overline{\mathbb{R}}$ be such that $m_1 \leq \partial^+ f(x)$ and $m_2 \leq \partial^- f(x)$.
    Then \eqref{ineq:estimate_of_sd-1} implies $\min\left\{ m_1, m_2 \right\} \leq m \leq f^{\sd}(x)$.
    \Cref{cor:estimate_of_sd-3} can be shown by applying \eqref{ineq:estimate_of_sd-1} similarly.
\end{proof}

The following examples illustrate the representation formula when one-sided derivatives may be infinite.

\begin{example} \label{ex:l.s.c._func}
    For fixed $\varepsilon \geq 0$, consider the function $f_{\varepsilon}:\mathbb{R} \to \mathbb{R}$ defined by 
    \begin{equation} \label{ex:l.s.c_func_epsilon}
        f_{\varepsilon}(x) =
        \begin{cases}
            x^2    &    \text{if } x \leq 1,    \\
            x + \varepsilon   &    \text{if } x > 1.
        \end{cases}
    \end{equation}
    Note that for any $\varepsilon \geq 0$, $f_{\varepsilon}'(1)$ does not exist, but $f_{\varepsilon}^{\sd}(1)$ exists.
    If $\varepsilon > 0$, then $f_{\varepsilon }$ is lower semicontinuous at $x = 1$ with $\partial^+ f_{\varepsilon}(1) = \infty$ and $\partial^- f_{\varepsilon}(1) = 2$,
    which implies that $f_{\varepsilon}^{\sd}(1) = \mathcal{C}(\infty, 2) = 2 + \sqrt{5} \approx 4.236 \ldots$.
    If $\varepsilon = 0$, then $\partial^+ f_{0}(1) = 1$ and $\partial^- f_{0}(1) = 2$, and hence $f_{0}^{\sd}(1) = \mathcal{C}(1, 2) \approx 1.387 \ldots$.
\end{example}

\begin{remark}  
    If one of the one-sided derivatives is infinite and the other is finite, then the specular derivative depends on the finite one and the sign of the infinite one, but not on how rapidly the difference quotient corresponding to the infinite one diverges.
    This is illustrated by $f_{\varepsilon}^{\sd}(1)$ for $\varepsilon > 0$, where $f_{\varepsilon}$ is defined in \eqref{ex:l.s.c_func_epsilon}.
\end{remark}

\begin{example} \label{ex:disc_func}
    We modify the value of $f_{\varepsilon}$ at $x=1$ as follows.
    For fixed $\varepsilon \geq 0$, consider $g_{\varepsilon}:\mathbb{R} \to \mathbb{R}$ defined by 
    \begin{displaymath} 
        g_{\varepsilon}(x) =
        \begin{cases}
            x^2    &    \text{if } x < 1,    \\
            \displaystyle 1 + \frac{1}{2} \varepsilon     &    \text{if } x = 1  , \\
            x + \varepsilon   &    \text{if } x > 1 .
        \end{cases}
    \end{displaymath}
    Note that $g_0 \equiv f_0$ on $\mathbb{R}$, and that $g_{\varepsilon}'(1)$ does not exist for any $\varepsilon \geq 0$.
    If $\varepsilon > 0$, then $g_{\varepsilon }$ is neither lower semicontinuous nor upper semicontinuous at $x = 1$.
    Moreover,
    \begin{displaymath}
        \partial^+ g_{\varepsilon}(1) 
        = \lim_{h \searrow 0} \left( 1 + \frac{\varepsilon }{2h} \right) 
        = \infty 
        = \lim_{h \searrow 0} \left( \frac{\varepsilon}{2h} + 2 - h \right)
        = \partial^- g_{\varepsilon}(1).
    \end{displaymath}
    Hence, $g_\varepsilon^{\sd}(1)$ does not exist. 
\end{example}

The examples in \cref{ex:l.s.c._func,ex:disc_func} show that changing the value of a function at a single point can affect both semicontinuity and specular differentiability.
For every $\varepsilon>0$, the functions $f_{\varepsilon}$ and $g_{\varepsilon}$ agree on $\mathbb{R}\setminus\{1\}$ but differ at $1$.
The function $f_{\varepsilon}$ is lower semicontinuous at $1$, whereas $g_{\varepsilon}$ is neither lower semicontinuous nor upper semicontinuous there.
Nevertheless, for every $h>0$,
\begin{displaymath}
    \frac{q^+_h f_{\varepsilon}(1) + q^-_h f_{\varepsilon}(1)}{2}
    =
    \frac{q^+_h g_{\varepsilon}(1) + q^-_h g_{\varepsilon}(1)}{2},
\end{displaymath}
and hence $f_{\varepsilon}^{\ast}(1)$ exists if and only if $g_{\varepsilon}^{\ast}(1)$ exists; in that case, $f_{\varepsilon}^{\ast}(1)=g_{\varepsilon}^{\ast}(1)$.
In contrast, $f_{\varepsilon}^{\sd}(1)$ exists, whereas $g_{\varepsilon}^{\sd}(1)$ does not.
This contrast illustrates the dependence of the specular derivative on the base-point value described in \cref{rmk:dependence_on_base_value}.

\begin{remark}
    The family in \cref{ex:l.s.c._func} also shows that uniform convergence of functions need not imply convergence of their
    specular derivatives.
    For each $n\in\mathbb{N}$, let $u_n:=f_{\frac{1}{n}}$.
    Then
    \begin{displaymath}
        \sup_{x\in\mathbb{R}}|u_n(x)-f_0(x)|
        =
        \frac{1}{n}
        \to
        0 
        \qquad\text{as}\qquad
        n \to \infty.
    \end{displaymath}
    Hence, $u_n\to f_0$ uniformly on $\mathbb{R}$.
    However, for every $n\in\mathbb{N}$,
    \begin{displaymath}
        u_n^{\sd}(1)
        =
        2 + \sqrt{5}
        \neq
        \mathcal{C}(1,2)
        =
        f_0^{\sd}(1).
\end{displaymath}
    Thus, $u_n^{\sd}(1)$ does not converge to $f_0^{\sd}(1)$.
\end{remark}

Finite one-sided differentiability provides the following sufficient condition for a function to belong to $S^1(I)$.

\begin{proposition} \label{prop:finite_one_sided}
    Let $I$ be an open interval in $\mathbb{R}$, and let $f:I\to\mathbb{R}$ be a function.
    If both $\partial^+f(x)$ and $\partial^-f(x)$ exist as finite real numbers for every $x\in I$, then $f\in S^1(I)$.
\end{proposition}

\begin{proof}
    Fix $x\in I$.
    By the continuity of $\mathcal{C}$ on $\mathbb{R}^2$, the limit in \eqref{eq:2nd_form_sd} exists and
    \begin{displaymath}
        f^{\sd}(x)
        =
        \mathcal{C}\left( \partial^+f(x), \partial^-f(x)  \right) \in \mathbb{R}.
    \end{displaymath}
    Since $x\in I$ was arbitrary, $f\in S^0(I)$.

    Moreover,
    \begin{displaymath}
        f(x+h)-f(x)
        =
        h \cdot q^+_h f(x) 
        \to 
        0 
    \end{displaymath}
    as $h\searrow0$.
    Hence, $f$ is right-continuous at $x$.
    Similarly, $f$ is left-continuous at $x$.
    Therefore, $f$ is continuous at $x$.

    Since $x\in I$ was arbitrary, $f\in C^0(I)$.
    Consequently, $f\in S^1(I)$.
\end{proof}

As an immediate consequence, every real-valued convex function $f$ on an open interval $I$ belongs to $S^1(I)$ and satisfies $f^{\sd}(x)\in\partial f(x)$ for every $x\in I$, where $\partial f(x)$ denotes the subdifferential of $f$ at $x$ in the sense of convex analysis.

\begin{corollary}
    \label{cor:convex_sd_subgradient}
    Let $I$ be an open interval in $\mathbb{R}$.
    If $f:I\to\mathbb{R}$ is convex, then $f\in S^1(I)$ and $f^{\sd}(x)\in\partial f(x)$ for every $x\in I$; equivalently,
    \begin{displaymath}
        f(y)
        \geq
        f(x)+f^{\sd}(x)(y-x)
    \end{displaymath}
    for all $x, y\in I$.
\end{corollary}

\begin{proof}
    Since $f$ is real-valued and convex on the open interval $I$, both $\partial^-f(x)$ and $\partial^+f(x)$ exist as finite real numbers for every $x\in I$; see \cite[Thm.~24.1]{1970_Rockafellar_BOOK}.
    Hence, \cref{prop:finite_one_sided} gives $f\in S^1(I)$.

    Fix $x\in I$.
    By \cref{cor:estimate_of_sd}~\ref{cor:estimate_of_sd-1} and the convexity of $f$,
    \begin{displaymath}
        \partial^-f(x)
        \leq
        f^{\sd}(x)
        \leq
        \partial^+f(x).
    \end{displaymath}
    Moreover, by the one-dimensional characterization of the subdifferential of a convex function \cite[Thm.~23.2]{1970_Rockafellar_BOOK}, we have $\partial f(x) = \left[ \partial^- f(x), \partial^+ f(x) \right]$.
    Therefore, $f^{\sd}(x)\in\partial f(x)$.
    Since $x\in I$ was arbitrary, the result follows.
\end{proof}

\section{Quasi-Mean Value Theorem}
\label{sec:Q-MVT}

The Quasi-Mean Value Theorem will serve as a key tool in establishing the regularity results in \cref{sec:regularity}.
Before proving the Quasi-Mean Value Theorem, we establish two preliminary results.
Throughout this section, all intervals are assumed to be nondegenerate; that is, whenever $[a,b]$ or $(a,b)$ appears, we assume that $a<b$.

\begin{lemma}
    \label{lem:existence_sd_one_sided_d}
    Let $f:[a, b] \to \mathbb{R}$ be continuous on $[a, b]$ and specularly differentiable on $(a, b)$.
    Then, the following statements hold:
    \begin{enumerate}[label=\upshape(\alph*)]
        \item \label{lem:existence_sd_one_sided_d-1} If $f(a) < f(b)$, then there exists $c_1 \in (a, b)$ such that $f^{\sd}(c_1)\geq 0$.
        Moreover, if at each $x\in(a,b)$ at least one of the one-sided derivatives $\partial^+f(x)$ and $\partial^-f(x)$ exists as an extended real number, then the same point $c_1$ also satisfies $\partial^+f(c_1) \geq 0$ and $\partial^-f(c_1) \geq 0$.
        \item \label{lem:existence_sd_one_sided_d-2} If $f(b) < f(a)$, then there exists $c_2 \in (a, b)$ such that $f^{\sd}(c_2) \leq 0$.
        Moreover, if at each $x\in(a,b)$ at least one of the one-sided derivatives $\partial^+f(x)$ and $\partial^-f(x)$ exists as an extended real number, then the same point $c_2$ also satisfies $\partial^+ f(c_2) \leq 0$ and $\partial^- f(c_2) \leq 0$.
    \end{enumerate}
\end{lemma}

\begin{proof}
    Without loss of generality, suppose that $f(a) < f(b)$; the case when $f(b) < f(a)$ can be shown similarly.
    Let $y \in \mathbb{R}$ be such that $f(a) < y < f(b)$.
    Define the set $X := \left\{x \in(a, b) \, \middle| \, y < f(x)\right\}$, which is bounded below by $a$; moreover, $a\notin X$.   
    The set $X$ is nonempty by the Intermediate Value Theorem.
    Since $X$ is bounded below and nonempty, there exists $c_1 \in [a, b)$ such that $c_1 = \inf X$.
    However, we have $c_1 \neq a$ by the continuity of $f$ and the fact that $y > f(a)$.
    Therefore, $c_1 \in (a, b)$.

    We claim that $f(c_1) = y$.
    Indeed, let $\varepsilon > 0$ be arbitrary.
    Since $f$ is continuous at $x = c_1$, there exists $\delta_1 > 0$ such that if $x \in (c_1 - \delta_1, c_1 + \delta_1)$, then $|f(c_1) - f(x)| < \varepsilon$.
    Since $(a, b)$ is open, there exists $\delta_2 > 0$ such that $(c_1 - \delta_2, c_1 + \delta_2) \subset (a, b)$.
    Choose $\delta_3 := \min \left\{ \delta_1, \delta_2 \right\}$.
    From $c_1 = \inf X$, there exists $x_1 \in X$ such that $c_1 \leq x_1 < c_1 + \delta_3$, and hence $f(c_1) > f(x_1) - \varepsilon  > y - \varepsilon$.
    For each $x_2 \in (c_1 - \delta_3, c_1)$, it follows that $x_2 \not \in X$ so that $f(c_1) < f(x_2) + \varepsilon \leq y + \varepsilon$.
    Combining these inequalities yields $y - \varepsilon < f(c_1) < y + \varepsilon$, which implies that $f(c_1) = y$.
    
    Since $f(c_1)=y$, we have $c_1\notin X$.
    Since $c_1 = \inf X$, there exists a sequence $h_n\searrow0$ such that $c_1+h_n\in X$ for every $n$.
    Hence,
    \begin{equation} \label{ineq:existence_right_d}
        \frac{f(c_1 + h_n) - f(c_1)}{h_n} > \frac{y - f(c_1)}{h_n} = 0
    \end{equation} 
    for all $n \in \mathbb{N}$.
    For sufficiently small $h>0$, we have 
    \begin{equation} \label{ineq:existence_left_d}
        \frac{f(c_1) - f(c_1 - h)}{h} \geq \frac{f(c_1) - y}{h} = 0
    \end{equation} 
    since $c_1 - h \notin X$.

    For sufficiently large $n$, take $h = h_n$ in \cref{ineq:existence_left_d}.    
    Combining this with \cref{ineq:existence_right_d} and applying \cref{lem:monotonicity_of_B}, we obtain
    \begin{displaymath}
        \mathcal{B}\left( \frac{f(c_1 + h_n) - f(c_1)}{h_n}, \frac{f(c_1) - f(c_1 - h_n)}{h_n} \right) 
        \geq
        \mathcal{B}(0, 0)
        =
        0
    \end{displaymath}
    for all sufficiently large $n$.
    Since $h_n\searrow0$ and $f$ is specularly differentiable at $c_1$, letting $n \to \infty$ in the preceding inequality and using the formula \eqref{eq:1st_form_sd} yields $f^{\sd}(c_1) \geq 0$.

    In addition, assume that at each $x\in(a,b)$ at least one of the one-sided derivatives $\partial^+f(x)$ and $\partial^-f(x)$ exists as an extended real number.
    Then both one-sided derivatives exist at every point of $(a,b)$ by \cref{lem:existence_one_sided_d-a,lem:existence_one_sided_d-b} of \cref{lem:existence_one_sided_d}.
    On the one hand, the inequality $\partial^+ f(c_1) \geq 0$ follows by taking the limit as $n \to \infty$ in \cref{ineq:existence_right_d}.
    On the other hand, the inequality $\partial^- f(c_1) \geq 0$ follows by taking the limit as $h \searrow 0$ in \cref{ineq:existence_left_d}.
\end{proof}

The second required statement for the proof of the Quasi-Mean Value Theorem is Quasi-Rolle's Theorem.
The proof of Quasi-Rolle's Theorem is inspired by \cite[Lem.~3]{1967_Aull} and \cite[Thm.~2.1.20]{2023_Jung} (see also \cite[Thm.~6.2]{1998_Sahoo_BOOK}).

\begin{theorem} 
    [Quasi-Rolle's Theorem] 
    \label{thm:quasi-Rolle}
    Let $f:[a, b] \to \mathbb{R}$ be continuous on $[a, b]$ and specularly differentiable on $(a, b)$.
    If $f(a) = f(b)$, then there exist $c_1, c_2 \in (a, b)$ such that
    \begin{displaymath} 
        f^{\sd}(c_1) \leq 0 \leq f^{\sd}(c_2).
    \end{displaymath}
\end{theorem}

\begin{proof}
    Define the function $g : [a, b] \to \mathbb{R}$ by 
    \begin{displaymath}
        g(x) := f(x) - f(a).
    \end{displaymath}
    Then, $g$ is continuous on $[a, b]$ and specularly differentiable on $(a, b)$ with $g(a) = 0 = g(b)$ and $g^{\sd} = f^{\sd}$ on $(a, b)$.
    If $g \equiv 0$ on $[a,b]$, then the conclusion is immediate.
    
    Now, suppose that $g \not\equiv 0$ on $[a, b]$.
    We consider the following three cases:
    \begin{enumerate}[label=\upshape(\roman*)]
        \item There exist $x_1, x_2 \in (a, b)$ such that $g(x_1) < 0 < g(x_2)$.
        \item $g \geq 0$ on $[a, b]$.
        \item $g \leq 0$ on $[a, b]$.
    \end{enumerate}
    If the first case holds, then $g(x_1) < g(a) = 0 < g(x_2)$.
    Applying \cref{lem:existence_sd_one_sided_d} to $g$ on $[a, x_1]$ and $[a, x_2]$, respectively, yields points $c_1 \in (a, x_1)$ and $c_2 \in (a, x_2)$ such that
    \begin{displaymath}
        g^{\sd}(c_1)\leq0\leq g^{\sd}(c_2).
    \end{displaymath}

    If the second case holds, then there exists $x_0 \in (a, b)$ such that $g(x_0) > 0$.
    Applying \cref{lem:existence_sd_one_sided_d} to $g$ on $[a,x_0]$ and $[x_0,b]$ yields the same conclusion.

    If the third case holds, then there exists $x_0 \in (a,b)$ such that $g(x_0)<0$.
    Applying \cref{lem:existence_sd_one_sided_d} to $g$ on $[a,x_0]$ and $[x_0,b]$ again yields the same conclusion.

    Since $g^{\sd}=f^{\sd}$ on $(a,b)$, the desired inequalities follow in all three cases.
\end{proof}

As in the classical proof of the Mean Value Theorem, Quasi-Rolle's Theorem yields the Quasi-Mean Value Theorem when $f(a)=f(b)$.
When $f(a)\neq f(b)$, however, subtracting the secant line does not reduce the problem to Quasi-Rolle's Theorem because specular differentiation is nonlinear.
We instead use the first and last crossings of suitable level sets, together with the monotonicity of $\mathcal{B}$ in each variable.

\begin{theorem}
    [Quasi-Mean Value Theorem] 
    \label{thm:q-MVT}
    Let $f:[a, b] \to \mathbb{R}$ be continuous on $[a, b]$ and specularly differentiable on $(a, b)$.
    Then, there exist $c_1, c_2 \in (a, b)$ such that
    \begin{equation} \label{ineq:q-MVT}
        f^{\sd}(c_1) \leq \frac{f(b) - f(a)}{b - a} \leq f^{\sd}(c_2).
    \end{equation}
\end{theorem}

\begin{proof}
    The case where $f(a) = f(b)$ follows from \cref{thm:quasi-Rolle}.
    Thus, we may assume that $f(a) \neq f(b)$.
    For convenience, write 
    \begin{displaymath}
        s := \frac{f(b) - f(a)}{b - a} \neq 0.
    \end{displaymath}
    Define the function $\varphi: [a, b] \to \mathbb{R}$ by 
    \begin{displaymath}
        \varphi(x) := f(x) - f(a) - s(x - a).
    \end{displaymath}
    Since $\varphi$ is continuous on $[a,b]$, it attains its minimum and maximum there; write
    \begin{displaymath}
        m := \min_{x\in[a,b]}\varphi(x) 
        \qquad\text{and}\qquad
        M := \max_{x\in[a,b]}\varphi(x).
    \end{displaymath}
    Then, we have $m \leq 0 \leq M$ since $\varphi(a) = \varphi(b) = 0$.
    We therefore distinguish the following three cases: $M > 0$, $M = 0 > m$, and $M = 0 = m$.
    
    First, suppose that $M > 0$.
    Let $\varepsilon \in \mathbb{R}$ be such that $0 < \varepsilon < M$.
    Define
    \begin{displaymath}
        X := \left\{ x \in [a, b] \, \middle| \, \varphi(x) > \varepsilon  \right\}.
    \end{displaymath}
    Then $X$ is nonempty, while $a, b \notin X$.
    Define $c_2 := \inf X$ and $c_1 := \sup X$.
    Then $a < c_2 \leq c_1 < b$, and the continuity of $\varphi$ gives $\varphi(c_1) = \varphi(c_2) = \varepsilon$.
    By the definition of $c_2$, there exists a sequence $h_n \searrow 0$ such that $c_2 + h_n\in X$.
    Moreover, $c_2 - h_n \notin X$ for all sufficiently large $n \in \mathbb{N}$.
    Hence
    \begin{displaymath}
        q^+_{h_n} \varphi(c_2) > 0 
        \qquad\text{and}\qquad
        q^-_{h_n} \varphi(c_2) \geq 0.
    \end{displaymath}
    Since $f(x)=\varphi(x)+f(a) + s(x-a)$, it follows that
    \begin{displaymath}
        q^+_{h_n} f(c_2) \geq s
        \qquad\text{and}\qquad
        q^-_{h_n} f(c_2) \geq s.
    \end{displaymath}
    By \cref{lem:monotonicity_of_B}, we have 
    \begin{displaymath}
        \mathcal{B}\left( q^+_{h_n} f(c_2), q^-_{h_n} f(c_2) \right)
        \geq
        \mathcal{B}(s, s)
        =
        s.
    \end{displaymath}
    Letting $n\to\infty$ yields $f^{\sd}(c_2)\geq s$.

    Similarly, there exists a sequence $k_n \searrow 0$ such that $c_1-k_n \in X$ and $c_1+k_n \notin X$ for all sufficiently large $n \in \mathbb{N}$.
    Therefore, $q^+_{k_n} \varphi(c_1) \leq 0$ and $q^-_{k_n} \varphi(c_1) < 0$.
    Repeating the preceding argument gives $f^{\sd}(c_1)\leq s$.
    Thus, the desired inequality \eqref{ineq:q-MVT} follows. 
    
    Next, suppose that $M = 0 > m$.
    Then, $\varphi$ is not identically zero. 
    Let $\varepsilon \in \mathbb{R}$ be such that $0 < \varepsilon < -m$, and define 
    \begin{displaymath}
        Y
        :=
        \left\{
            x \in [a, b] \;\middle|\; \varphi(x) < -\varepsilon
        \right\}.
    \end{displaymath}
    Since $Y$ is nonempty and $\varphi(a) = \varphi(b) = 0 > -\varepsilon$, we have $a, b \notin Y$.
    Define $c_1:=\inf Y$ and $c_2:=\sup Y$.
    Then $a < c_1 \leq c_2 < b$.
    Moreover, the continuity of $\varphi$ gives $\varphi(c_1) = \varphi(c_2) = -\varepsilon$.
    Choose $h_n \searrow 0$ such that $c_1 + h_n \in Y$ and $c_1 - h_n \notin Y$.
    Then
    \begin{displaymath}
        q^+_{h_n} f(c_1) < s
        \qquad\text{and}\qquad
        q^-_{h_n} f(c_1) \leq s.
    \end{displaymath}
    Hence, by \cref{lem:monotonicity_of_B},
    \begin{displaymath}
        \mathcal{B}\left( q^+_{h_n} f(c_1), q^-_{h_n} f(c_1) \right)
        \leq
        \mathcal{B}(s, s)
        =
        s.
    \end{displaymath}
    Letting $n\to\infty$ gives $f^{\sd}(c_1)\leq s$.

    Similarly, choose $k_n\searrow0$ such that
    $c_2-k_n\in Y$ and $c_2+k_n\notin Y$.
    Then
    \begin{displaymath}
        q^+_{k_n} f(c_2) \geq s
        \qquad\text{and}\qquad
        q^-_{k_n} f(c_2) > s.
    \end{displaymath}
    Therefore,
    \begin{displaymath}
        \mathcal{B}\left( q^+_{k_n} f(c_2), q^-_{k_n} f(c_2) \right)
        \geq
        \mathcal{B}(s, s)
        =
        s.
    \end{displaymath}
    Letting $n\to\infty$ gives $s \leq f^{\sd}(c_2)$.
    Thus, \eqref{ineq:q-MVT} follows.

    Finally, if $M = 0 = m$, then $\varphi \equiv 0$.
    Thus, $f(x)=f(a) + s(x-a)$ on $[a, b]$, and hence $f^{\sd}(x) = s$ for every $x \in (a, b)$.
    Consequently, any points $c_1, c_2 \in (a, b)$ satisfy the desired inequality.
\end{proof}

\begin{remark}
    Quasi-Rolle's Theorem and the Quasi-Mean Value Theorem do not require the existence of one-sided derivatives.
\end{remark}

Monotonicity can be characterized by the sign of the specular derivative as follows.

\begin{corollary}
    \label{cor:monotonicity_and_sd}
    Let $f:[a, b] \to \mathbb{R}$ be continuous on $[a, b]$ and specularly differentiable on $(a, b)$.
    \begin{enumerate}[label=\upshape(\alph*)]
        \item \label{cor:monotonicity_and_sd-1} $f^{\sd}(x) \geq 0$ (resp. $f^{\sd}(x) \leq 0$) for all $x \in (a, b)$ if and only if $f$ is nondecreasing (resp. nonincreasing) on $[a, b]$.
        \item \label{cor:monotonicity_and_sd-2} $f^{\sd}(x) = 0$ for all $x \in (a, b)$ if and only if $f$ is constant on $[a, b]$.
    \end{enumerate}
\end{corollary}

\begin{proof}
    We first prove \cref{cor:monotonicity_and_sd-1}.
    Suppose that $f^{\sd}(x) \geq 0$ for all $x \in (a, b)$.
    Let $x_1, x_2 \in [a, b]$ with $x_1 < x_2$.
    By \cref{thm:q-MVT} applied to $f$ on $[x_1, x_2]$, there exists $c_1 \in (x_1, x_2)$ such that
    \begin{displaymath}
        f^{\sd}(c_1)(x_2 - x_1) \leq f(x_2) - f(x_1).
    \end{displaymath}
    Since $f^{\sd}(c_1) \geq 0$ and $x_2 - x_1>0$, we obtain $f(x_2) \geq f(x_1)$.
    Thus, $f$ is nondecreasing on $[a,b]$.

    Conversely, suppose that $f$ is nondecreasing on $[a,b]$.
    Fix $x \in (a,b)$.
    For every sufficiently small $h>0$, we have $f(x+h)-f(x) \geq 0$ and $f(x)-f(x-h) \geq 0$.
    Then, \cref{lem:A}~\ref{lem:A-7} implies 
    \begin{displaymath}
        \mathcal{C} \left( q^+_h f(x), q^-_h f(x) \right) \geq 0.
    \end{displaymath}
    Taking the limit as $h \searrow 0$, we obtain $f^{\sd}(x)\geq 0$ from the formula \eqref{eq:2nd_form_sd}.
    Therefore, $f^{\sd}(x) \geq 0$ for all $x \in(a, b)$.

    The case $f^{\sd}(x) \leq 0$ is proved similarly, using the right-hand inequality in \cref{thm:q-MVT} for the sufficient part and the monotonicity of a nonincreasing function for the necessary part.

    Finally, \cref{cor:monotonicity_and_sd-2} follows from \cref{cor:monotonicity_and_sd-1}.
    If $f^{\sd}(x)=0$ for all $x\in(a, b)$, then $f^{\sd}(x) \geq 0$ and $f^{\sd}(x) \leq 0$ for all $x \in (a, b)$.
    By \cref{cor:monotonicity_and_sd-1}, $f$ is both nondecreasing and nonincreasing on $[a, b]$, and hence $f$ is constant.
    Conversely, if $f$ is constant on $[a, b]$, then $q^+_h f(x) = q^-_h f(x) = 0$ for every $x \in (a, b)$ and every sufficiently small $h > 0$.
    Therefore, $f^{\sd}(x) = 0$ for all $x \in (a, b)$ by \eqref{eq:2nd_form_sd} and \cref{lem:A-1}.
\end{proof}

\section{Discontinuities of specularly differentiable functions}
\label{sec:discontinuities}

For comparison, Charzyński \cite[Thm.~1]{1933_Charzynski} proved that if a real-valued function on an open interval has a finite symmetric derivative at every point, then its set of discontinuities is scattered and hence at most countable.
Freiling \cite[Thm.~6]{1990_Freiling} later extended this conclusion by allowing countably many exceptional points at which the function is only assumed to be symmetrically continuous.
No analogous conclusion holds for everywhere specularly differentiable functions: \cref{ex:nowhere_continuous} satisfies $\Disc(f)=\mathbb{R}$ and $f^{\sd}\equiv0$.
The next result nevertheless shows that the discontinuities at which the specular derivative is nonzero form an at most countable set.

For the proof of the countability result below, we introduce notation for points at which a function is strictly increasing or decreasing across the base point.
For a real-valued function $f$ defined on an open interval $I$ in $\mathbb{R}$, define
\begin{align*}
    \Inc(f)
    &:=
    \left\{
        x\in I
        \,\middle|\,
        \begin{array}{l}
            \text{there exists }r>0\text{ such that } (x-r,x+r)\subset I \text{ and } \\
            f(y)<f(x)<f(z) \text{ for all }y\in(x-r,x) \text{ and }z\in(x,x+r)
        \end{array}
    \right\},\\
    \Dec(f)
    &:=
    \Inc(-f).
\end{align*}

We now state the main result of this section.

\begin{theorem} \label{thm:countable_nonzero_sd_discontinuities}
    Let $I$ be an open interval in $\mathbb{R}$.
    If $f \in S^0(I)$, then the set 
    \begin{displaymath}
        \left\{ x \in \Disc(f) \, \middle| \, f^{\sd}(x) \neq 0 \right\}
    \end{displaymath}
    is at most countable.
\end{theorem}

\begin{proof}
    For $\varepsilon > 0$, define the sets 
    \begin{displaymath}
        E^+_{\varepsilon} := \left\{ x \in I \, \middle| \, f^{\sd}(x) \geq \varepsilon  \right\} 
        \qquad\text{and}\qquad
        E^-_{\varepsilon} := \left\{ x \in I \, \middle| \, f^{\sd}(x) \leq -\varepsilon  \right\},
    \end{displaymath}
    and the functions $\varphi^+_{\varepsilon}, \varphi^-_{\varepsilon} : I \to \mathbb{R}$ by 
    \begin{displaymath}
        \varphi^+_{\varepsilon}(x) := f(x) + \frac{x}{\varepsilon} 
        \qquad\text{and}\qquad
        \varphi^-_{\varepsilon}(x) := f(x) - \frac{x}{\varepsilon} .
    \end{displaymath}
    
    Fix $\varepsilon > 0$.
    We claim that $E^+_{\varepsilon} \subset \Inc(\varphi^+_{\varepsilon})$ and $E^-_{\varepsilon} \subset \Dec(\varphi^-_{\varepsilon})$.
    Observe that, for $x \in I$, 
    \begin{equation}    \label{eq:countable_nonzero_sd_discontinuities-1}
        \lim_{h \searrow 0} \left( \arctan \bigl( q_h^+ f(x) \bigr)+\arctan \bigl( q_h^- f(x) \bigr) \right) 
        =
        2\arctan\bigl(f^{\sd}(x)\bigr)
    \end{equation}
    from the definition of $\mathcal{B}$ and the formula \eqref{eq:1st_form_sd}.
    
    First, let $x \in E^+_{\varepsilon}$.
    Since
    \begin{displaymath}
        2\arctan\bigl(f^{\sd}(x)\bigr)
        \geq
        2\arctan\varepsilon
        >
        \arctan\varepsilon,
    \end{displaymath}
    we have, for all sufficiently small $h>0$,
    \begin{displaymath}
        \arctan \bigl( q_h^+ f(x) \bigr)+\arctan \bigl( q_h^- f(x) \bigr)
        >
        \arctan\varepsilon
    \end{displaymath}
    by \eqref{eq:countable_nonzero_sd_discontinuities-1}.
    Since $\arctan t<\frac{\pi}{2}$ for every $t\in\mathbb{R}$, we obtain
    \begin{displaymath}
        \arctan \bigl( q_h^{\pm} f(x) \bigr)
        >
        \arctan\varepsilon-\frac{\pi}{2}
        =
        -\arctan \frac{1}{\varepsilon} ,
    \end{displaymath}
    and hence $q_h^{\pm} f(x)>-\frac{1}{\varepsilon}$.
    For all sufficiently small $h>0$,
    \begin{align*}
        \varphi^+_{\varepsilon}(x+h) - \varphi^+_{\varepsilon}(x)
        &=
        h\left(q_h^+ f(x)+\frac{1}{\varepsilon}\right)>0,\\
        \varphi^+_{\varepsilon}(x) - \varphi^+_{\varepsilon}(x-h)
        &=
        h\left(q_h^- f(x)+\frac{1}{\varepsilon}\right)>0.
    \end{align*}
    Therefore, $x\in\Inc(\varphi^+_{\varepsilon})$.
    This proves $E^+_{\varepsilon} \subset \Inc(\varphi^+_{\varepsilon})$.

    Second, let $x \in E^-_{\varepsilon}$.
    Then, for all sufficiently small $h>0$,
    \begin{displaymath}
        \arctan \bigl( q_h^+ f(x) \bigr)+\arctan \bigl( q_h^- f(x) \bigr)
        <
        -\arctan\varepsilon,
    \end{displaymath}
    which can be proved from \eqref{eq:countable_nonzero_sd_discontinuities-1} as before.
    Since $\arctan t > -\frac{\pi}{2}$ for every $t\in\mathbb{R}$, we obtain
    \begin{displaymath}
        \arctan \bigl( q_h^{\pm} f(x) \bigr)
        <
        -\arctan\varepsilon+\frac{\pi}{2}
        =
        \arctan \frac{1}{\varepsilon} ,
    \end{displaymath}
    and hence $q_h^{\pm} f(x)<\frac{1}{\varepsilon}$.
    For all sufficiently small $h>0$,
    \begin{align*}
        \varphi^-_{\varepsilon}(x+h) - \varphi^-_{\varepsilon}(x)
        &=
        h\left(q_h^+ f(x)-\frac{1}{\varepsilon}\right) < 0,\\
        \varphi^-_{\varepsilon}(x) - \varphi^-_{\varepsilon}(x-h)
        &=
        h\left(q_h^- f(x) - \frac{1}{\varepsilon}\right) < 0.
    \end{align*}
    Therefore, $x\in\Dec(\varphi^-_{\varepsilon})$.
    This shows $E^-_{\varepsilon} \subset \Dec(\varphi^-_{\varepsilon})$.

    Next, note that $\Disc(\varphi^+_{\varepsilon}) = \Disc(f) = \Disc(\varphi^-_{\varepsilon})$.
    Consequently,
    \begin{align*}
        \Disc(f)\cap E^+_{\varepsilon}
        &\subset
        \Disc(\varphi^+_{\varepsilon})
        \cap
        \Inc(\varphi^+_{\varepsilon})
        \subset
        \Disc(\varphi^+_{\varepsilon})
        \cap
        \bigl(
            \Inc(\varphi^+_{\varepsilon})
            \cup
            \Dec(\varphi^+_{\varepsilon})
        \bigr),\\
        \Disc(f)\cap E^-_{\varepsilon}
        &\subset
        \Disc(\varphi^-_{\varepsilon})
        \cap
        \Dec(\varphi^-_{\varepsilon})
        \subset
        \Disc(\varphi^-_{\varepsilon})
        \cap
        \bigl(
            \Inc(\varphi^-_{\varepsilon})
            \cup
            \Dec(\varphi^-_{\varepsilon})
        \bigr).
    \end{align*}
    By \cref{lem:countable_discontinuous_strict_increase}, both sets on the left-hand sides are at most countable.
    Since
    \begin{displaymath}
        \left\{ x \in \Disc(f) \,\middle|\, f^{\sd}(x)\neq 0 \right\}
        =
        \bigcup_{n=1}^{\infty} \left( \left( \Disc(f) \cap E^+_{\frac{1}{n}} \right) \cup \left( \Disc(f) \cap E^-_{\frac{1}{n}} \right) \right),
    \end{displaymath}
    the right-hand side is a countable union of at most countable sets, and is therefore at most countable.
\end{proof}

The restriction to points at which the specular derivative is nonzero is sharp: the following example satisfies $\Disc(f)=\mathbb{R}$ and $f^{\sd}\equiv0$.
Earlier Hamel-basis examples from the theory of symmetric differentiation provide the starting point for the construction.
Uher \cite[p.~37]{19871988_Uher} gave an example exhibiting extreme behavior under symmetric differentiation on a set of full outer measure, and Freiling and Rinne \cite[p.~519]{19881989_Freiling} later presented several related nonmeasurable examples.
In one of their constructions, they considered $g=\chi_G$, where $G$ is the kernel of a Hamel-coordinate functional; the symmetric derivative of $g$ is zero on $G$ and does not exist at any point of $\mathbb{R}\setminus G$.
The following example modifies this construction.

Define the function $v:\mathbb{Q}\to\mathbb{Z}\cup\{\infty\}$ by
\begin{equation}    \label{def:v}
    v(q)
    :=
    \begin{cases}
        k
        &
        \begin{aligned}[t]
            &\text{if $q\neq0$, where $k\in\mathbb{Z}$ is the unique integer for which}\\
            &\text{there exist odd integers $a$ and $b$ satisfying $qb=2^k a$} ,
        \end{aligned}\\
        \infty
        & \text{if $q=0$}.
\end{cases}
\end{equation}
Define the function $\rho:\mathbb{Q}\to[0,\infty)$ by
\begin{equation}    \label{def:rho}
    \rho(q)
    :=
    \begin{cases}
        2^{-v(q)}
        & \text{if $q\neq0$},\\
        0
        & \text{if $q=0$}.
    \end{cases}
\end{equation}
The functions $v$ and $\rho$ are the $2$-adic valuation and the $2$-adic absolute value on $\mathbb{Q}$, respectively; see \cite[Chap.~2]{2020_Gouvea_BOOK} for the general definitions and their basic properties.
We record below the properties of $v$ and $\rho$ needed in the sequel.
Note that, for all $q, r \in \mathbb{Q} \setminus \{0\}$,
\begin{equation}    \label{eq:v}
    v(qr) = v(q) + v(r).
\end{equation}
Note that, for all $q, r\in\mathbb{Q}$,
\begin{equation}    \label{eq:rho}
    \rho(qr) = \rho(q)\rho(r)
    \qquad\text{and}\qquad
    \rho(q+r) \leq \max\{\rho(q), \rho(r)\}.
\end{equation}
Moreover, if $\rho(q)\neq\rho(r)$, then
\begin{equation}    \label{eq:rho_strict}
    \rho(q+r)
    =
    \rho(q-r)
    =
    \max\{\rho(q), \rho(r)\}.
\end{equation}

We next construct the map used in the following example.
For a set $E\subset\mathbb{R}$, let $\Span_{\mathbb{Q}}(E)$ denote the set of all finite linear combinations of elements of $E$ with coefficients in $\mathbb{Q}$.
Let $\mathcal{H}$ be a Hamel basis of $\mathbb{R}$ over $\mathbb{Q}$, and fix an element $e_0 \in \mathcal{H}$.
Define the function $T:\mathbb{R}\to\mathbb{Q}$ by
\begin{equation}    \label{def:T}
    T(x) := q,
\end{equation}
where $q \in \mathbb{Q}$ and $y \in\Span_{\mathbb{Q}}\bigl(\mathcal{H}\setminus\{e_0\}\bigr)$ are the unique elements satisfying $x = q e_0 + y$.
The map $T$ is linear in the following sense:
\begin{equation}    \label{eq:T_linearity}
    T(ax_1+bx_2)
    =
    aT(x_1)+bT(x_2)
\end{equation}
for every $x_1,x_2\in\mathbb{R}$ and $a,b\in\mathbb{Q}$.
Moreover, $T$ is surjective since $T(qe_0)=q$ for every $q\in\mathbb{Q}$.
Also, $\mathcal{H}\setminus\{e_0\}\neq\varnothing$, since otherwise $\mathbb{R} = \Span_{\mathbb{Q}}(\{e_0\}) = \left\{ qe_0 \, \middle| \, q\in\mathbb{Q} \right\}$ would be countable.
Therefore, $\ker T = \Span_{\mathbb{Q}} \bigl(\mathcal{H}\setminus\{e_0\}\bigr) \neq \{0\}$, and hence $\ker T$ is nontrivial.

\begin{example}[A specularly differentiable function that is nowhere continuous] \label{ex:nowhere_continuous}
    Define the function $f:\mathbb{R}\to\mathbb{R}$ by
    \begin{displaymath}
        f(x) := -\rho(T(x)),
    \end{displaymath}
    where $\rho$ and $T$ are defined by \eqref{def:rho} and \eqref{def:T}, respectively.
    We claim that $f$ is specularly differentiable on $\mathbb{R}$ and $f^{\sd}(x) = 0$ for every $x\in\mathbb{R}$.
    Nevertheless, $f$ is nowhere continuous and is not Lebesgue measurable.

    First, fix $x\in\mathbb{R}$.
    We show that $f^{\sd}(x)=0$.
    For $h>0$, set $\alpha_h := q_h^+f(x)$ and $\beta_h := q_h^-f(x)$.
    By \eqref{eq:T_linearity}, we have 
    \begin{equation}
        \label{eq:2-adic-difference-quotients}
        \alpha_h
        =
        \frac{\rho(T(x))-\rho(T(x) + T(h))}{h} 
        \quad\text{and}\quad
        \beta_h
        =
        \frac{\rho(T(x) - T(h))-\rho(T(x))}{h}
    \end{equation}
    for every $h > 0$.
    We claim that 
    \begin{equation} \label{ex:nowhere_continuous-2}
        \lim_{h \searrow 0} \left( \arctan\alpha_h+\arctan\beta_h \right)
        =
        0.
    \end{equation}

    If $x\in\ker T$, then $\rho(-T(h))=\rho(T(h))$, and hence \eqref{eq:2-adic-difference-quotients} gives $\beta_h = -\alpha_h$ for every $h>0$.
    Since the arctangent function is odd, \eqref{ex:nowhere_continuous-2} holds.
    
    Now suppose that $x \notin \ker T$, that is, $T(x) \neq 0$.
    Let $k:=v(T(x))$, and write
    \begin{displaymath}
        T(x) = 2^k \cdot \frac{a}{b},
    \end{displaymath}
    where $a$ and $b$ are odd integers.
    For each $h>0$, either $\rho(T(x))\neq\rho(T(h))$ or $\rho(T(x))=\rho(T(h))$.
    On the one hand, suppose that $\rho(T(x))\neq\rho(T(h))$.
    Then
    \begin{displaymath}
        \rho(T(x) + T(h))
        =
        \rho(T(x) - T(h))
        =
        \max\{\rho(T(x)), \rho(T(h))\}
    \end{displaymath}
    by \eqref{eq:rho_strict}.
    It follows from \eqref{eq:2-adic-difference-quotients} that $\beta_h=-\alpha_h$, and hence
    \begin{displaymath}
        \arctan\alpha_h+\arctan\beta_h=0.
    \end{displaymath}

    On the other hand, suppose that $\rho(T(x)) = \rho(T(h))$.
    Since $x \notin \ker T$, we have $T(x) \neq 0$, and hence $\rho(T(x)) > 0$.
    The equality $\rho(T(x)) = \rho(T(h))$ therefore implies that $T(h) \neq 0$.
    Since $\rho(s) = 2^{-v(s)} > 0$ for every $s \in \mathbb{Q} \setminus \left\{ 0 \right\}$, we obtain $k = v(T(h))$.
    Write 
    \begin{displaymath}
        T(h)=2^k \cdot \frac{c}{d},
    \end{displaymath}
    where $c$ and $d$ are odd integers.
    For either choice of sign, if $ad\pm bc=0$, then $\rho\bigl(T(x)\pm T(h)\bigr)=0$.
    Otherwise, since $ad\pm bc$ is even and $bd$ is odd, \eqref{eq:v} gives
    \begin{align*}
        v\bigl(T(x)\pm T(h)\bigr)
        &= v\left( 2^k \cdot \frac{ad \pm bc}{bd} \right) 
        = v \bigl(2^k \bigr) + v \left( ad \pm bc \right) + v \left( \frac{1}{bd} \right)  \\
        &= k + v \left( ad \pm bc \right) 
        \geq k + 1.
    \end{align*}
    Therefore, in either case,
    \begin{equation}  \label{ex:nowhere_continuous-1}
        \rho\bigl(T(x) \pm T(h)\bigr) 
        \leq
        2^{-(k+1)}
        = \frac{\rho(T(x))}{2}.
    \end{equation}
    Hence, \eqref{eq:2-adic-difference-quotients} and \eqref{ex:nowhere_continuous-1} yield
    \begin{displaymath}
        \alpha_h
        \geq
        \frac{\rho(T(x))}{2h}
        >
        0
        \qquad\text{and}\qquad
        \beta_h
        \leq
        -\frac{\rho(T(x))}{2h}
        <
        0.
    \end{displaymath}
    Therefore,
    \begin{align*}
        \left|
            \arctan\alpha_h+\arctan\beta_h
        \right|
        &\leq
        \left(\frac{\pi}{2}-\arctan\alpha_h\right) + \left(\frac{\pi}{2}+\arctan\beta_h\right) \\
        &=
        \arctan\frac{1}{\alpha_h} + \arctan\frac{1}{|\beta_h|}
        \leq
        \frac{1}{\alpha_h} + \frac{1}{|\beta_h|} \\
        &\leq 
        \frac{4h}{\rho(T(x))}.  
    \end{align*}

    Combining the two possibilities, we have
    \begin{displaymath}
        \left|
            \arctan\alpha_h+\arctan\beta_h
        \right|
        \leq
        \frac{4h}{\rho(T(x))}
    \end{displaymath}
    for every $h>0$: the left-hand side is zero when $\rho(T(x))\neq\rho(T(h))$, while the preceding estimate applies when $\rho(T(x))=\rho(T(h))$.
    Sending $h \searrow 0$ proves \eqref{ex:nowhere_continuous-2}.

    Therefore, using \cref{lem:ABC} and \eqref{ex:nowhere_continuous-2}, we obtain
    \begin{displaymath}
        f^{\sd}(x)
        =
        \lim_{h \searrow 0} \mathcal{C}(\alpha_h,\beta_h)
        =
        \lim_{h \searrow 0} \tan\left( \frac{ \arctan\alpha_h + \arctan\beta_h} {2} \right)
        =
        0.
    \end{displaymath}
    Since $x\in\mathbb{R}$ was arbitrary, $f^{\sd} \equiv 0$ on $\mathbb{R}$.

    We next prove that $f$ is nowhere continuous.
    Since $\ker T$ is nontrivial, choose $u \in \ker T\setminus\{0\}$.
    By \eqref{eq:T_linearity}, $T(qu)=qT(u)=0$ for every $q\in\mathbb{Q}$, and hence $\mathbb{Q}u\subset\ker T$.
    Moreover, $\mathbb{Q}u$ is dense in $\mathbb{R}$ since $u \neq 0$ and $\mathbb{Q}$ is dense in $\mathbb{R}$.
    Hence, $\ker T$ is dense in $\mathbb{R}$.

    We also have
    \begin{displaymath}
        T^{-1}(\{0\})=\ker T
        \qquad\text{and}\qquad
        T^{-1}(\{1\})=e_0+\ker T.
    \end{displaymath}
    Since a translate of a dense set is dense, both
    $T^{-1}(\{0\})$ and $T^{-1}(\{1\})$ are dense in $\mathbb{R}$.
    Furthermore, since $\rho(0)=0$ and $\rho(1)=1$,
    \begin{displaymath}
        f=0 \quad\text{on }T^{-1}(\{0\})
        \qquad\text{and}\qquad
        f=-1 \quad\text{on }T^{-1}(\{1\}).
    \end{displaymath}
    Thus, continuity at any $x\in\mathbb{R}$ would imply, by approaching
    $x$ through these two dense sets, that $f(x)=0$ and $f(x)=-1$,
    a contradiction.
    Therefore, $f$ is nowhere continuous.

    Finally, we prove that $f$ is not Lebesgue measurable.
    Since $\{0\}$ is a Borel set and $f^{-1}(\left\{ 0 \right\}) = \ker T$, it is enough to show that $\ker T$ is not Lebesgue measurable.

    Suppose, to the contrary, that $\ker T$ were Lebesgue measurable.
    First, assume that $\mathcal{L}^1(\ker T) = 0$.
    For every $x\in\mathbb{R}$, \eqref{eq:T_linearity} and $T(e_0)=1$ give
    \begin{displaymath}
        T\bigl(x-T(x)e_0\bigr)
        =
        T(x)-T(x)T(e_0)
        =
        0.
    \end{displaymath}
    Therefore, $x - T(x) e_0 \in \ker T$.
    Since $T(x)\in\mathbb{Q}$, it follows that
    \begin{displaymath}
        \mathbb{R}
        =
        \bigcup_{q\in\mathbb{Q}} \bigl(qe_0+\ker T\bigr).
    \end{displaymath}
    By translation invariance, every set $qe_0+\ker T$ is null.
    Since $\mathbb{Q}$ is countable, this would imply that $\mathcal{L}^1(\mathbb{R}) = 0$, a contradiction.
    
    Now assume that $\mathcal{L}^1(\ker T)>0$.
    By Steinhaus' theorem (see \cite{2020_Sadhukhan}), the set $\left\{ u-v \, \middle| \, u,v\in\ker T \right\}$ contains a neighborhood of $0$.
    Since $\ker T$ is an additive subgroup of $\mathbb{R}$, we have $u - v\in\ker T$ whenever $u, v\in\ker T$.
    Therefore, $\ker T$ itself contains a neighborhood of $0$.
    That is, there exists $\delta>0$ such that $(-\delta, \delta)\subset\ker T$.
    Given $x\in\mathbb{R}$, choose $n\in\mathbb{N}$ sufficiently large that $\left|\frac{x}{n}\right|<\delta$.
    Then $\frac{x}{n}\in\ker T$, and hence
    \begin{displaymath}
        x
        =
        \sum_{i=1}^{n}\frac{x}{n}
        \in\ker T,
    \end{displaymath}
    since $\ker T$ is an additive subgroup.
    Thus, $\ker T=\mathbb{R}$, contradicting $T(e_0)=1$.

    Hence, $\ker T$ is not Lebesgue measurable, and therefore neither is $f$.
\end{example}

\section{Regularity of specularly differentiable functions}
\label{sec:regularity}

Uher established several regularity results for symmetrically differentiable functions.
If the symmetric derivative $f^{\ast}$ exists almost everywhere on a Lebesgue measurable set $E\subset\mathbb{R}$, then $f$ is classically differentiable almost everywhere on $E$, and hence $f'=f^{\ast}$ almost everywhere on $E$; see \cite[Cor.~1]{19821983_Uher}.
Moreover, the existence of $f^{\ast}$ almost everywhere on $\mathbb{R}$ already implies that $f$ is Lebesgue measurable;
see \cite[Cor.~2]{19821983_Uher}.

Uher's conclusion cannot, in general, be strengthened to classical differentiability outside a countable set.
Foran \cite{1977_Foran} constructed a continuous, everywhere symmetrically differentiable function whose classical derivative fails to exist on an uncountable set.
Ponomarev \cite{1985_Ponomarev} subsequently strengthened this phenomenon by constructing a continuous, everywhere symmetrically differentiable function whose set of points of classical nondifferentiability is precisely a nonempty perfect set $P$, while the function is of class $C^1(\mathbb{R}\setminus P)$.

In contrast, the existence of the specular derivative, even everywhere, does not in general imply either Lebesgue measurability or classical differentiability almost everywhere.
Motivated by this contrast, we investigate below additional assumptions that ensure regularity of specularly differentiable functions and their specular derivatives.

Recalling \cref{thm:countable_nonzero_sd_discontinuities}, we obtain the following sufficient condition for $f$ to be continuous almost everywhere.
Without this additional assumption, the conclusion may fail; see \cref{ex:nowhere_continuous}.

\begin{theorem}
    Let $I$ be an open interval in $\mathbb{R}$.
    If $f \in S^0 (I)$ and $f^{\sd}(x) \neq 0$ for $\mathcal{L}^1$-a.e. $x \in I$, then $f$ is continuous $\mathcal{L}^1$-a.e. on $I$.
\end{theorem}

\begin{proof}
    Writing 
    \begin{displaymath}
        E := \left\{ x\in\Disc(f) \, \middle| \, f^{\sd}(x)\neq 0 \right\} 
        \qquad\text{and}\qquad
        F := \left\{ x\in\Disc(f) \, \middle| \, f^{\sd}(x)=0 \right\},
    \end{displaymath}
    we have $\Disc(f) = E \cup F$.
    By \cref{thm:countable_nonzero_sd_discontinuities}, we have $\mathcal{L}^1(E) = 0$.
    Moreover, $\mathcal{L}^1(F) = 0$ since $F$ is contained in the $\mathcal{L}^1$-null set $\left\{ x\in I \,\middle|\, f^{\sd}(x)=0 \right\}$.
    Therefore
    \begin{displaymath}
        0 \leq \mathcal{L}^1 \bigl(\Disc(f)\bigr) \leq \mathcal{L}^1(E) + \mathcal{L}^1(F) = 0,
    \end{displaymath}
    which implies that $\mathcal{L}^1 \bigl(\Disc(f)\bigr) = 0$.
\end{proof}

We next study the regularity of the specular derivative, beginning with measurability.
If a specularly differentiable function is Lebesgue measurable, then so is its specular derivative.
A substantially stronger conclusion is available in the symmetric setting.
Khintchine \cite[p.~217]{1927_Khintchine} proved that if $f:\mathbb{R}\to\mathbb{R}$ is Lebesgue measurable, then $f'(x)$ exists as a finite number at almost every point $x$ satisfying
\begin{displaymath}
    \limsup_{h \to 0}
    \frac{f(x+h)-f(x-h)}{2h}
    <
    +\infty.
\end{displaymath}

\begin{proposition}
    Let $I$ be an open interval in $\mathbb{R}$.
    If $f \in S^0(I)$ and $f$ is Lebesgue measurable, then $f^{\sd}$ is Lebesgue measurable.
\end{proposition}

\begin{proof}
    Let $\widetilde{f}:\mathbb{R}\to\mathbb{R}$ be the extension of $f$ by zero outside $I$.
    Since $I$ is open and $f$ is Lebesgue measurable on $I$, the function $\widetilde{f}$ is Lebesgue measurable on $\mathbb{R}$.
    For every $n\in\mathbb{N}$, define the function $g_n:I\to\mathbb{R}$ by
    \begin{displaymath}
        g_n(x)
        :=
        \mathcal{C}\left(
            \frac{\widetilde{f}\left(x+\frac{1}{n}\right)-\widetilde{f}(x)}{\frac{1}{n}},
            \frac{\widetilde{f}(x)-\widetilde{f}\left(x-\frac{1}{n}\right)}{\frac{1}{n}}
        \right).
    \end{displaymath}
    Since translations preserve Lebesgue measurability and $\mathcal{C}$ is continuous on $\mathbb{R}^2$, each $g_n$ is Lebesgue measurable on $I$.

    Fix $x\in I$.
    Since $I$ is open, we have $x - \frac{1}{n}, x+\frac{1}{n} \in I$ for all sufficiently large $n\in\mathbb{N}$, and hence
    \begin{displaymath}
        g_n(x)
        =
        \mathcal{C}\left(
            \frac{f\big(x+\frac{1}{n}\big)-f(x)}{\frac{1}{n}},
            \frac{f(x)-f\big(x-\frac{1}{n}\big)}{\frac{1}{n}}
        \right).
    \end{displaymath}
    Therefore, the specular differentiability of $f$ at $x$ implies that $g_n(x) \to f^{\sd}(x)$ as $n \to \infty$.
    Hence, $f^{\sd}$ is the pointwise limit of a sequence of Lebesgue measurable functions and is therefore Lebesgue measurable on $I$.
\end{proof}

Under the stronger assumption that $f$ is continuous, the preceding measurability conclusion can be improved.
In fact, the specular derivative of every function in $S^1(I)$ is of Baire class $1$.
For symmetric derivatives, Larson \cite[Thm.~2.1]{1983_Larson} proved that if a function is symmetrically differentiable everywhere, then its symmetric derivative is of Baire class $1$.

\begin{proposition} \label{prop:Baire_one}
    Let $I$ be an open interval in $\mathbb{R}$.
    If $f \in S^1(I)$, then $f^{\sd} \in B^1(I)$.
\end{proposition}

\begin{proof}
    Choose a continuous function $\varphi : I \to (0, \infty)$ such that $[x - \varphi(x), x + \varphi(x)]\subset I$ for every $x \in I$.
    For each $n \in \mathbb{N}$, define the functions $h_n : I \to (0, \infty)$ and $g_n:I\to\mathbb{R}$ by 
    \begin{displaymath}
        h_n(x) := \frac{\varphi(x)}{n} 
    \end{displaymath}
    and 
    \begin{displaymath}
        g_n(x)
        :=
        \mathcal{C}\left( \frac{f(x+h_n(x))-f(x)}{h_n(x)}, \frac{f(x)-f(x-h_n(x))}{h_n(x)} \right),
    \end{displaymath}
    respectively.
    Since $f$, $h_n$, and $\mathcal{C}$ are continuous and $h_n$ is positive, each $g_n$ is continuous on $I$.
    For every fixed $x\in I$, we have $h_n(x)\searrow0$, and hence the formula \eqref{eq:2nd_form_sd} yields
    \begin{displaymath}
        \lim_{n\to\infty}g_n(x)=f^{\sd}(x).
    \end{displaymath}
    Thus, $f^{\sd}$ is the pointwise limit of a sequence of continuous real-valued functions on $I$, and therefore $f^{\sd}\in B^1(I)$.
\end{proof}

The following lemma shows that a function in $S^1(I)$ is classically differentiable at every continuity point of its specular derivative.

\begin{lemma}
    \label{lem:differentiability_at_continuity_point_of_sd}
    Let $I$ be an open interval in $\mathbb{R}$, and let $f\in S^1(I)$.
    If $f^{\sd}$ is continuous at $x\in I$, then $f$ is classically differentiable at $x$ and $f'(x)=f^{\sd}(x)$.
\end{lemma}

\begin{proof}
    Fix $x \in I$, and choose $r > 0$ such that $[x - r, x + r] \subset I$.
    Let $\varepsilon > 0$.
    Since $f^{\sd}$ is continuous at $x$, there exists $\delta \in (0, r)$ such that
    \begin{displaymath}
        \left\vert f^{\sd}(x)-f^{\sd}(y) \right\vert <\varepsilon
    \end{displaymath}
    whenever $|x - y| < \delta$.
    For $h \in (0, \delta)$, the Quasi-Mean Value Theorem applied on $[x,x+h]$ yields points $c_1(h), c_2(h) \in (x, x+h)$ such that
    \begin{displaymath}
        f^{\sd}(c_1(h))
        \leq
        q^+_h f(x)
        \leq
        f^{\sd}(c_2(h)).
    \end{displaymath}
    Since $|c_i(h)-x|<h<\delta$ for $i=1,2$, it follows that
    \begin{displaymath}
        f^{\sd}(x)-\varepsilon
        <
        q^+_h f(x)
        <
        f^{\sd}(x)+\varepsilon.
    \end{displaymath}
    Hence, $\partial^+f(x)=f^{\sd}(x)$.

    Similarly, applying \cref{thm:q-MVT} on $[x-h,x]$ gives points
    $c_3(h),c_4(h)\in(x-h,x)$ such that
    \begin{displaymath}
        f^{\sd}(c_3(h))
        \leq
        q^-_h f(x)
        \leq
        f^{\sd}(c_4(h)).
    \end{displaymath}
    Since $|c_i(h)-x|<h<\delta$ for $i=3,4$, we obtain
    \begin{displaymath}
        f^{\sd}(x)-\varepsilon
        <
        q^-_h f(x)
        <
        f^{\sd}(x)+\varepsilon.
    \end{displaymath}
    Therefore, $\partial^-f(x)=f^{\sd}(x)$.
    Thus, $\partial^+f(x)=\partial^-f(x)=f^{\sd}(x)$, which proves that $f'(x)$ exists and equals $f^{\sd}(x)$.
\end{proof}

Combining the Baire class $1$ regularity of the specular derivative with the preceding lemma yields residual classical differentiability.
For symmetric derivatives, Belna, Evans, and Humke \cite[Thm.~4]{1978_Belna} proved the stronger result that every measurable function $f:\mathbb{R}\to\mathbb{R}$ that is symmetrically differentiable everywhere is classically differentiable outside a $\sigma$-porous set, which is necessarily meager.

\begin{theorem}
    \label{thm:residual_differentiability}
    Let $I$ be an open interval in $\mathbb{R}$.
    If $f\in S^1(I)$, then $\Cont(f^{\sd})$ is a dense $G_\delta$ subset of $I$.
    Moreover, $f$ is classically differentiable at every point of $\Cont(f^{\sd})$, and $f'=f^{\sd}$ on $\Cont(f^{\sd})$.
\end{theorem}

\begin{proof}
    By \cref{prop:Baire_one}, $f^{\sd}\in B^1(I)$.
    By the classical theorem on points of continuity of Baire class $1$ functions \cite[Thm.~24.14]{1995_Kechris}, $\Cont(f^{\sd})$ is a comeager $G_\delta$ subset of $I$.
    Since $I$ is locally compact and Hausdorff, it is a Baire space by the Baire category theorem \cite[Thm.~8.4]{1995_Kechris}.
    Since every comeager subset of a Baire space is dense \cite[Prop.~8.1]{1995_Kechris}, $\Cont(f^{\sd})$ is dense in $I$.
    By \cref{lem:differentiability_at_continuity_point_of_sd}, the function $f$ is classically differentiable at every point of $\Cont(f^{\sd})$ and $f'=f^{\sd}$ on $\Cont(f^{\sd})$.
\end{proof}

\begin{remark}
    The exceptional set $I\setminus\Cont(f^{\sd})$ in \cref{thm:residual_differentiability} need not be countable.
    Indeed, let $f$ be the continuous function constructed by Ponomarev \cite{1985_Ponomarev}, and let $P$ be the nonempty perfect set such that $f$ is not classically differentiable at any point of $P$, $f^{\ast}(x) = 0$ for every $x\in P$, and $f\in C^1(\mathbb{R}\setminus P)$.
    In particular, $f$ is symmetrically differentiable on $\mathbb{R}$.

    We first show that $f\in S^1(\mathbb{R})$.
    If $x\in P$, then \cref{lem:A}~\ref{lem:A-9} gives
    \begin{displaymath}
        \left| \mathcal{C}\left( q^+_h f(x), q^-_h f(x) \right) \right|
        \leq 
        \left\vert \frac{q^+_h f(x) + q^-_h f(x)}{2} \right\vert .
    \end{displaymath}
    The right-hand side converges to zero as $h \searrow 0$, since $f^{\ast}(x)=0$.
    Hence $f^{\sd}(x)=0$.
    If $x\notin P$, then $f$ is classically differentiable at $x$, and therefore $f^{\sd}(x)=f'(x)$.
    Thus $f$ is specularly differentiable on $\mathbb{R}$.
    Since $f$ is continuous, it follows that $f\in S^1(\mathbb{R})$.

    Moreover, we have $\Disc(f^{\sd})=P$.
    Indeed, since $P$ is closed, $\mathbb{R}\setminus P$ is open.
    On this open set, $f^{\sd}=f'$, and $f'$ is continuous.
    Hence $f^{\sd}$ is continuous at every point of $\mathbb{R}\setminus P$.
    On the other hand, if $f^{\sd}$ were continuous at some $x\in P$, then \cref{lem:differentiability_at_continuity_point_of_sd} would imply that $f$ is classically differentiable at $x$, contradicting the construction.
    Therefore, $f^{\sd}$ is discontinuous at every point of $P$, proving the equality.

    Since every nonempty perfect subset of $\mathbb{R}$ is uncountable, the exceptional set $\mathbb{R}\setminus\Cont(f^{\sd})=P$ is uncountable.
\end{remark}

We next establish an almost-everywhere differentiability result under the assumption that the specular derivative is bounded.
A continuous function with a bounded specular derivative is Lipschitz continuous and hence classically differentiable $\mathcal{L}^1$-almost everywhere.
An analogous statement for symmetric derivatives can be found in \cite[Thm.~4]{1967_Aull}.

\begin{proposition}
    \label{prop:cont_bounded_sd_impy_Lip_cont}
    Let $I$ be an open interval in $\mathbb{R}$.
    If $f \in S^1(I)$ and $f^{\sd}$ is bounded on $I$, then $f$ is Lipschitz continuous on $I$, $f'$ exists $\mathcal{L}^1$-a.e. on $I$, and $f' = f^{\sd}$ $\mathcal{L}^1$-a.e. on $I$.
\end{proposition}

\begin{proof}
    Since $f^{\sd}$ is bounded on $I$, there exists $M>0$ such that $|f^{\sd}(x)|\leq M$ for every $x \in I$.
    Let $x_1,x_2\in I$ be such that $x_1 < x_2$.
    By \cref{thm:q-MVT}, there exist $c_1, c_2 \in (x_1, x_2)$ such that
    \begin{displaymath}
        -M 
        \leq
        f^{\sd}(c_1)
        \leq
        \frac{f(x_2) - f(x_1)}{x_2 - x_1}
        \leq
        f^{\sd}(c_2)
        \leq
        M,
    \end{displaymath}
    which implies $|f(x_1) - f(x_2)| \leq M |x_1 - x_2|$.
    Since $x_1$ and $x_2$ were arbitrary, $f$ is Lipschitz continuous on $I$ with Lipschitz constant $M$.

    By Rademacher's theorem (see, for example, \cite[Thm.~3.2]{2015_Evans}), $f$ is differentiable $\mathcal{L}^1$-a.e. on $I$.
    At every point at which $f'$ exists, the specular derivative agrees with the classical derivative.
    Therefore, the desired result follows.
\end{proof}

Finally, we turn to the case of a continuous specular derivative.
In this setting, specular differentiability can be upgraded to classical differentiability everywhere.
An analogous result for symmetric derivatives was proved by Aull; see \cite[Thm.~2]{1967_Aull}.

\begin{theorem}
    \label{thm:cont_sd_implies_d}
    Let $I$ be an open interval in $\mathbb{R}$.
    If $f\in S^1(I)$ and $f^{\sd}\in C^0(I)$, then $f \in C^1(I)$ and $f'=f^{\sd}$ on $I$.
\end{theorem}

\begin{proof}
    By \cref{lem:differentiability_at_continuity_point_of_sd}, we have $f'=f^{\sd}$ on $I$.
    Since $f^{\sd}\in C^0(I)$, it follows that $f'\in C^0(I)$, and hence $f\in C^1(I)$.
\end{proof}

\section{Second-order specular differentiability}
\label{sec:2nd_order}

In this section, we introduce the second-order specular derivative and investigate its relation to classical regularity.
We show that even a continuous function may be twice specularly differentiable without being classically differentiable. 
We then identify additional continuity conditions under which classical regularity is recovered and compare the resulting function classes with $C^1(I)$ and $C^2(I)$.

\begin{definition}  \label{def:S2}
    Let $I$ be an open interval in $\mathbb{R}$, and let $f\in S^0(I)$.
    If $f^{\sd}$ is specularly differentiable at $x\in I$, then the
    \emph{second-order specular derivative} of $f$ at $x$ is defined by
    \begin{displaymath}
        f^{\sd\sd}(x):=(f^{\sd})^{\sd}(x).
    \end{displaymath}
    We say that $f$ is \emph{twice specularly differentiable} on $I$
    if $f^{\sd}\in S^0(I)$.
    We define the class 
    \begin{displaymath}
        S^2 (I) := \left\{ f \in S^1(I) \, \middle| \, f^{\sd} \in S^1(I) \right\}.
    \end{displaymath}    
\end{definition}

We emphasize that twice specular differentiability alone does not imply that $f\in S^2(I)$.
Indeed, $f\in S^2(I)$ precisely when $f$ is twice specularly differentiable and both $f$ and $f^{\sd}$ are continuous.

\begin{example}[A function in $C^1$ but not in $S^2$] \label{ex:C1notS2}
    Consider the function $f:\mathbb{R}\to\mathbb{R}$ defined by
    \begin{displaymath}
        f(x):=
        \begin{cases}
            \displaystyle x^3\sin\frac{1}{x}
            & \text{if }x\neq0, \\[2mm]
            0
            & \text{if }x=0.
        \end{cases}
    \end{displaymath}
    Then $f\in C^1(\mathbb{R})$ and $f^{\sd}=f'$.

    However, for $h > 0$, let $A_h := 3h\sin\frac{1}{h} - \cos\frac{1}{h}$.
    Since $f^{\sd}=f'$, the limit defining the second-order specular derivative at $0$ is
    \begin{displaymath}
        \lim_{h \searrow 0} \mathcal{C}\left( q^+_h (f') (0), q^-_h (f') (0) \right)
        =
        \lim_{h \searrow 0} \mathcal{C}\left( A_h, A_h \right)
        =
        \lim_{h \searrow 0} A_h
    \end{displaymath}
    which does not exist.
    Thus, $f^{\sd\sd}(0)$ does not exist, and therefore $f\notin S^2(\mathbb{R})$.  
\end{example}

\begin{example}[A function in $S^2$ but not in $C^2$] \label{ex:S2notC2}
    Consider the function $f:\mathbb{R}\to\mathbb{R}$ defined by $f(x) := \frac{1}{2}x|x|$.
    A direct calculation gives $f^{\sd}(x)=|x|$ and $f^{\sd\sd}(x)=\sgn(x)$.
    Hence, $f\in S^2(\mathbb{R})$.
    However, $f''(0)$ does not exist, and therefore $f \notin C^2(\mathbb{R})$.
\end{example}

The following theorem summarizes the regularity results for specularly differentiable functions.

\begin{theorem} \label{thm:chain_regularity}
    Let $I$ be a nonempty open interval in $\mathbb{R}$.
    Then,
    \begin{displaymath}
        C^2(I) \subsetneq S^2(I) \subsetneq C^1(I) \subsetneq S^1(I) \subsetneq C^0(I).
    \end{displaymath}    
\end{theorem}

\begin{proof}
    The inclusion $S^2(I)\subset C^1(I)$ follows from \cref{thm:cont_sd_implies_d}.
    Applying \cref{prop:finite_one_sided} to a function in $C^1(I)$ gives $C^1(I)\subset S^1(I)$.
    If $f\in C^2(I)$, then applying the same proposition to $f$ and $f'$, together with the identity $f^{\sd}=f'$, gives $f\in S^2(I)$.
    The inclusion $S^1(I)\subset C^0(I)$ follows directly from the definition.
    To prove strictness for an arbitrary nonempty open interval $I$, fix $x_0\in I$, translate each example in \cref{ex:S1notC1}, \cref{ex:C0notS1}, \cref{ex:C1notS2}, and \cref{ex:S2notC2} so that its exceptional point is $x_0$, and then restrict it to $I$.
\end{proof}

Applying \cref{prop:finite_one_sided} to the first-order specular derivative yields the following second-order consequence.

\begin{proposition}
    \label{prop:finite_2nd_one_sided}
    Let $I$ be an open interval in $\mathbb{R}$, and let $f\in S^0(I)$.
    Suppose that both $\partial^+(f^{\sd})(x)$ and $\partial^-(f^{\sd})(x)$ exist as finite real numbers for every $x \in I$.
    Then $f^{\sd}\in S^1(I)$.
    In particular, $f$ is twice specularly differentiable on $I$.

    If, in addition, $f\in S^1(I)$, then $f\in S^2(I)$.
    In this case, $f\in C^1(I)$ and $f'=f^{\sd}$ on $I$.
\end{proposition}

\begin{proof}
    Applying \cref{prop:finite_one_sided} to $f^{\sd}$ gives $f^{\sd}\in S^1(I)$.
    Hence, $f$ is twice specularly differentiable on $I$.

    If $f\in S^1(I)$, then the definition of $S^2(I)$ gives $f\in S^2(I)$.
    The remaining conclusions follow from \cref{thm:cont_sd_implies_d}.
\end{proof}

The condition $f^{\sd}\in C^0(I)$ in \cref{thm:cont_sd_implies_d} cannot be replaced by the mere existence of $f^{\sd\sd}$ on $I$.
Although \cref{prop:finite_2nd_one_sided} provides a sufficient condition for this continuity, its finiteness assumption is not
automatic, even for a continuous twice specularly differentiable function.
Indeed, the following example is continuous and twice specularly differentiable, but its specular derivative is discontinuous and has neither one-sided derivative at the origin, while the function itself is not classically differentiable there.

\begin{example}[A twice specularly differentiable function that is continuous but not classically differentiable]   \label{ex:notS2notC1}
    Consider the function $f : \mathbb{R} \to \mathbb{R}$ defined by
    \begin{displaymath}
        f(x)
        :=
        \begin{cases}
            \displaystyle x^2\sin\frac{1}{x} + \frac{|x|}{4} - \frac{x^2}{8}\sin\frac{2}{|x|} 
            & \text{if }x \neq 0, \\[2mm]
            0
            & \text{if }x = 0.
        \end{cases}
    \end{displaymath}
    We claim that $f$ is twice specularly differentiable on $\mathbb{R}$, but $f\notin S^2(\mathbb{R})$ and $f\notin C^1(\mathbb{R})$.
    Since $f \in C^0(\mathbb{R})$ and $f \in C^{\infty}(\mathbb{R} \setminus \left\{ 0 \right\})$, the second-order specular derivative $f^{\sd\sd}(x)$ exists for every $x \neq 0$.
    Consequently, it is enough to show that $f^{\sd\sd}(0)$ exists, whereas $f'(0)$ does not exist and $f^{\sd}$ is not continuous at $0$.

    First, $f'(0)$ does not exist since $\partial^+f(0)=\frac{1}{4}$ and $\partial^-f(0)=-\frac{1}{4}$.
    On the other hand, by \eqref{eq:2nd_form_sd} and the continuity of $\mathcal{C}$ on $\mathbb{R}^2$, $f^{\sd}(0) = \mathcal{C}\left(\frac{1}{4},-\frac{1}{4}\right) = 0$.
    
    Second, since $f' = f^{\sd}$ on $\mathbb{R}\setminus\{0\}$, 
    \begin{displaymath}
        f^{\sd}(x)
        =
        \begin{cases}
            \displaystyle 2x\sin\frac{1}{x} -\cos\frac{1}{x} -\frac{x}{4}\sin\frac{2}{|x|} + \frac{\sgn(x)}{4}\left( 1 + \cos \frac{2}{|x|} \right)
            & \text{if }x\neq0, \\[2mm]
            0
            & \text{if }x=0.
        \end{cases}
    \end{displaymath}
    In particular, for every $n\geq1$, $f^{\sd}\left(\frac{1}{2 \pi n}\right) = -\frac{1}{2}$ and $f^{\sd}\left(\frac{1}{(2n + 1)\pi}\right) = \frac{3}{2}$.
    Thus, $f^{\sd}$ is discontinuous at $0$.

    Third, neither $\partial^+(f^{\sd})(0)$ nor $\partial^-(f^{\sd})(0)$ exists, even as an extended real number.
    Indeed, for each $n \in \mathbb{N}$, let $h_n := \frac{1}{2\pi n}$ and $k_n := \frac{1}{(2n+1)\pi}$.
    Then $f^{\sd}(h_n)=-\frac{1}{2}$, $f^{\sd}(k_n)=\frac{3}{2}$, $f^{\sd}(-h_n)=-\frac{3}{2}$, and $f^{\sd}(-k_n)=\frac{1}{2}$.
    Consequently, as $n \to \infty$, 
    \begin{displaymath}
        q^+_{h_n} (f^{\sd}) (0) = -\pi n \to -\infty 
        \qquad\text{and}\qquad
        q^+_{k_n} (f^{\sd}) (0) = \frac{3(2n+1)\pi}{2} \to \infty,
    \end{displaymath}
    whereas
    \begin{displaymath}
        q^-_{h_n}(f^{\sd})(0) = 3\pi n \to \infty
        \qquad\text{and}\qquad
        q^-_{k_n}(f^{\sd})(0) = -\frac{(2n+1)\pi}{2} \to -\infty.
    \end{displaymath}
    Thus, neither one-sided derivative exists.

    Finally, it remains to show that $f^{\sd\sd}(0)$ nevertheless exists.
    From \eqref{eq:2nd_form_sd}, 
    \begin{displaymath}
        f^{\sd \sd}(0)
        = \lim_{h \to 0} \mathcal{C}\left( \frac{f^{\sd}(h)}{h}, -\frac{f^{\sd}(-h)}{h} \right).
    \end{displaymath}
    For $h \neq 0$, let 
    \begin{displaymath}
        A_h := \frac{f^{\sd}(h)}{h}
        \qquad\text{and}\qquad
        B_h := -\frac{f^{\sd}(-h)}{h}.
    \end{displaymath}
    Since $A_{-h}=B_h$, $B_{-h}=A_h$, and $\mathcal{C}$ is symmetric, it suffices to consider $h \searrow 0$.
    For $h > 0$, let 
    \begin{displaymath}
        u_h := \frac{1}{h}\cos \frac{1}{h} - 2\sin\frac{1}{h}
        \qquad\text{and}\qquad
        v_h:= \frac{1}{2}\cos \frac{1}{h} \left( \frac{1}{h} \cos \frac{1}{h} - \sin \frac{1}{h} \right).
    \end{displaymath}
    Then, $A_h = v_h - u_h$, $B_h = v_h + u_h$, and $ v_h = \frac{1}{2} \cos \frac{1}{h} \left(u_h+\sin\frac{1}{h}\right)$.

    Let $\varepsilon>0$ be arbitrary.
    Since
    \begin{displaymath}
        \lim_{R \to \infty} \tan\left[
        \frac{1}{2}\left(
            \frac{\pi}{2}
            -
            \arctan\left(\frac{R-1}{2}\right)
        \right)
    \right] 
    = 0,
    \end{displaymath}
    choose $R>1$ sufficiently large that
    \begin{displaymath}
        \tan\left[ \frac{1}{2}\left( \frac{\pi}{2} - \arctan\left(\frac{R - 1}{2}\right) \right) \right]
        <
        \varepsilon.
    \end{displaymath}
    Next, choose $\delta>0$ such that $\delta < \frac{2\varepsilon}{(R+1)(R+2)}$, and let $0 < h < \delta$.
    On the one hand, if $|u_h|\leq R$, then
    \begin{displaymath}
        |v_h|
        \leq
        \frac{1}{2} \left\vert \cos \frac{1}{h} \right\vert  \left\vert u_h + \sin\frac{1}{h} \right\vert
        =
        \frac{h}{2} \left\vert u_h + 2 \sin \frac{1}{h} \right\vert \left\vert u_h + \sin\frac{1}{h} \right\vert
        \leq 
        \frac{h(R+1)(R+2)}{2}.
    \end{displaymath}
    Therefore, by \cref{lem:A}~\ref{lem:A-9},
    \begin{displaymath}
        |\mathcal{C}(A_h,B_h)|
        \leq
        \frac{|A_h+B_h|}{2}
        =
        |v_h|
        \leq
        \frac{h(R+1)(R+2)}{2}
        <
        \varepsilon.
    \end{displaymath}
    On the other hand, suppose now that $|u_h|>R$.
    Since
    \begin{displaymath}
        |v_h|
        \leq 
        \frac{1}{2} \left\vert \cos \frac{1}{h} \right\vert  \left( |u_h| + \left\vert \sin \frac{1}{h} \right\vert  \right)
        \leq
        \frac{|u_h|+1}{2}
        <
        |u_h|
    \end{displaymath}
    and $A_h B_h = v_h^2 - u_h^2 < 0$, the numbers $A_h$ and $B_h$ have opposite signs.
    Writing $m_h := \min\{|A_h|, |B_h|\}$ and $M_h := \max\{|A_h|, |B_h|\}$, we have 
    \begin{displaymath}
        m_h
        \geq
        |u_h|-|v_h|
        \geq
        \frac{|u_h|-1}{2}
        >
        \frac{R-1}{2}.
    \end{displaymath}
    Consequently,
    \begin{align*}
        \left| \arctan A_h + \arctan B_h \right|
        &= \left| \arctan |A_h| - \arctan |B_h| \right| 
        = \arctan M_h - \arctan m_h  \\
        &\leq \frac{\pi}{2} - \arctan m_h 
        \leq \frac{\pi}{2} - \arctan \left( \frac{R - 1}{2} \right).
    \end{align*}    
    Using $\mathcal{B} = \mathcal{C}$, we obtain
    \begin{displaymath}
        |\mathcal{C}(A_h, B_h)|
        \leq
        \tan\left[ \frac{1}{2} \left( \frac{\pi}{2} - \arctan\left(\frac{R-1}{2}\right) \right) \right]
        <\varepsilon.
    \end{displaymath}
    Thus, in either case, $|\mathcal{C}(A_h,B_h)| < \varepsilon$ whenever $0 < h < \delta$.
    Therefore, $\mathcal{C}(A_h, B_h)$ converges to zero as $h \searrow 0$.
    Hence, $f^{\sd\sd}(0) = 0$.

    Thus, twice specular differentiability in the present sense does not ensure continuity of $f^{\sd}$ or classical differentiability of $f$.
\end{example}

\appendix 
\crefalias{section}{appendix}

\section{Analysis of auxiliary functions}
\label{apx:A_and_B}

This appendix presents further analysis of the auxiliary functions $\mathcal{A}$, $\mathcal{B}$, and $\mathcal{C}$ defined in \cref{def:sd,def:BC}.
We first give the deferred proof of \cref{lem:ABC}.

\begin{proof}[Proof of \cref{lem:ABC}]
    First, let $(a, b, c) \in \mathbb{R} \times \mathbb{R} \times (0, \infty)$ be fixed.
    We claim that $\mathcal{A}(a, b, c) = \mathcal{B} \left( \frac{a}{c}, \frac{b}{c} \right)$.
    Consider the substitution $a = c \tan \theta_1 $ and $b = c \tan \theta_2$ for some $\theta_1, \theta_2 \in \left( - \frac{\pi}{2}, \frac{\pi}{2} \right)$.
    Then, we find that  
    \begin{align*}
        \mathcal{A}(a, b, c)
        &= \frac{c^2 \tan \theta_1 \sqrt{\tan^2 \theta_2 + 1} + c^2 \tan \theta_2 \sqrt{\tan^2 \theta_1 + 1}}{c^2 \sqrt{\tan^2 \theta_2 + 1} + c^2 \sqrt{\tan^2 \theta_1 + 1}}  \\
        &= \frac{\tan \theta_1 \sec \theta_2 + \tan \theta_2 \sec \theta_1}{\sec \theta_2 + \sec \theta_1} 
        = \frac{\sin \theta_1 + \sin \theta_2}{\cos \theta_1 + \cos \theta_2} \\
        &= \frac{2 \sin \left( \frac{\theta_1 + \theta_2}{2} \right) \cos \left( \frac{\theta_1 - \theta_2}{2} \right)}{2 \cos \left( \frac{\theta_1 + \theta_2}{2} \right) \cos \left( \frac{\theta_1 - \theta_2}{2} \right)}
        = \tan \left( \frac{\theta_1 + \theta_2}{2} \right) 
        = \mathcal{B} \left( \frac{a}{c}, \frac{b}{c} \right).
    \end{align*}
    Here, we have used the facts that $\cos\theta_1>0$, $\cos\theta_2>0$, and $\cos\left(\frac{\theta_1-\theta_2}{2}\right) > 0$.

    Next, let $(\alpha, \beta) \in \overline{\mathbb{R}}^2$ be fixed. 
    We claim that $\mathcal{B}(\alpha, \beta) = \mathcal{C}(\alpha, \beta)$.
    If $(\alpha, \beta) \in \mathbb{R}^2$, then 
    \begin{align*}
        \mathcal{B}(\alpha, \beta)
        &= \frac{\sin \left( \frac{1}{2} \arctan \alpha + \frac{1}{2} \arctan \beta \right)\cos \left( \frac{1}{2}\arctan \alpha - \frac{1}{2} \arctan \beta   \right)}{\cos \left( \frac{1}{2} \arctan \alpha + \frac{1}{2} \arctan \beta  \right)\cos \left( \frac{1}{2} \arctan \alpha - \frac{1}{2}\arctan \beta  \right)}   \\
        &= \frac{\sin(\arctan \alpha) + \sin(\arctan \beta )}{\cos(\arctan \alpha) + \cos(\arctan \beta )}   \\
        &= \left(\frac{\alpha}{\sqrt{1 + \alpha^2}} + \frac{\beta}{\sqrt{1 + \beta^2}}\right)\left(\frac{1}{\sqrt{1 + \alpha^2}} + \frac{1}{\sqrt{1 + \beta^2}}\right)^{-1} 
        = \mathcal{C}(\alpha, \beta).
    \end{align*}
    
    Suppose that exactly one of $\alpha$ and $\beta$ is infinite.
    For every $x\in\mathbb{R}$, the half-angle formula for the tangent function yields the identity
    \begin{equation}\label{eq:repr_of_sd}
        \tan\left(
            \frac{1}{2}\arctan x \pm \frac{\pi}{4}
        \right)
        =
        x\pm\sqrt{1+x^2}.
    \end{equation}
    Thus, for every $x \in \mathbb{R}$, $\mathcal{B}(x,\pm\infty) = x\pm\sqrt{1+x^2} = \mathcal{C}(x,\pm\infty)$.
    By symmetry, the same conclusion holds for $(\alpha,\beta)=(\pm\infty,x)$.

    Finally, the remaining cases, in which both $\alpha$ and $\beta$ are infinite, follow directly from the definitions.
\end{proof}

We next establish further properties of the auxiliary functions.
First, the function $\mathcal{B}$ is strictly increasing in each of its arguments on $\overline{\mathbb{R}}^2$.

\begin{lemma}   \label{lem:monotonicity_of_B}
    If $(\alpha_1, \beta_1), (\alpha_2, \beta_2) \in \overline{\mathbb{R}}^2$ with $\alpha_1 \leq \alpha_2$ and $\beta_1 \leq \beta_2$, then 
    \begin{displaymath} 
        \mathcal{B}(\alpha_1,\beta_1)
        \leq
        \mathcal{B}(\alpha_2,\beta_2).
    \end{displaymath}
    Moreover, the inequality is strict if at least one of the two inequalities is strict.
\end{lemma}

\begin{proof}
    By the conventions in \eqref{def:B}, the maps $\arctan:\overline{\mathbb{R}} \to \big[-\frac{\pi}{2},\frac{\pi}{2}\big]$ and $\tan : \big[-\frac{\pi}{2},\frac{\pi}{2}\big] \to\overline{\mathbb{R}}$ are strictly increasing.
    Hence, the assumptions imply
    \begin{displaymath}
        \frac{\arctan \alpha_1 +\arctan \beta_1}{2}
        \leq
        \frac{\arctan \alpha_2 +\arctan \beta_2}{2}.
    \end{displaymath}
    Applying the extended tangent function proves the desired inequality.

    If $\alpha_1 < \alpha_2$ or $\beta_1 < \beta_2$, then the corresponding inequality between the averaged angles, and hence the
    resulting inequality, is strict.
\end{proof}

Next, the function $\mathcal{C}$ belongs to $C^{\infty}(\mathbb{R}^2)$.
Indeed, the identities
\begin{displaymath}
    \frac{x}{\sqrt{1+x^2}}
    =
    \sin(\arctan x)
    \qquad\text{and}\qquad
    \frac{1}{\sqrt{1+x^2}}
    =
    \cos(\arctan x)
\end{displaymath}
show that the numerator and denominator in the definition of $\mathcal{C}$ are smooth on $\mathbb{R}^2$.
Moreover, $\cos(\arctan\alpha)+\cos(\arctan\beta) > 0$ for every $(\alpha,\beta)\in\mathbb{R}^2$.
Thus, the denominator never vanishes, and hence $\mathcal{C}\in C^{\infty}(\mathbb{R}^2)$.

Observe that $\mathcal{C}$ admits the following alternative representation on $\mathbb{R}^2$:
\begin{equation}    \label{eq:algebraic_repr_C}   
    \mathcal{C}(\alpha, \beta) 
    =     
    \begin{cases}
        \displaystyle \frac{\alpha \beta - 1 + \sqrt{(1 + \alpha^2)(1 + \beta^2)}}{\alpha + \beta}    
        &    \text{if } \alpha + \beta \neq 0,    \\[0.5em]
        0    
        &    \text{if } \alpha + \beta = 0 .
    \end{cases}
\end{equation}

We collect several useful properties of the function $\mathcal{C}$.

\begin{lemma}   \label{lem:A}
    Let $(\alpha, \beta) \in \mathbb{R}^2$ be arbitrary. 
    Then the following properties hold.
    \begin{enumerate}[label=\upshape(\alph*)]
        \item $\mathcal{C}(\alpha, \alpha) = \alpha$. \label{lem:A-1}
        \item $\mathcal{C}(\alpha, \beta) = \mathcal{C}(\beta, \alpha)$.    \label{lem:A-2}
        \item $\mathcal{C}(-\alpha, -\beta) = -\mathcal{C}(\alpha, \beta)$. \label{lem:A-3}
        \item $\beta \leq \mathcal{C}(\alpha, \beta) \leq \alpha$ if $\beta \leq \alpha$. \label{lem:A-4}
        \item $\mathcal{C}(\alpha, \beta) = \alpha$ or $\mathcal{C}(\alpha, \beta) = \beta$ if and only if $\alpha = \beta$. \label{lem:A-5}
        \item $\mathcal{C}(\alpha, \beta) \neq 0$ if and only if $\alpha + \beta \neq 0$.    \label{lem:A-6}
        \item $\sgn(\alpha + \beta) = \sgn(\mathcal{C}(\alpha, \beta))$.  \vspace{1mm}  \label{lem:A-7}
        \item $\displaystyle \frac{\alpha + \beta - |\alpha + \beta|}{4} \leq \mathcal{C}(\alpha, \beta) \leq \frac{\alpha + \beta + |\alpha + \beta|}{4}$. \vspace{1mm} \label{lem:A-8}
        \item $\displaystyle -\frac{|\alpha + \beta|}{2}\leq \mathcal{C}(\alpha, \beta) \leq \frac{|\alpha + \beta|}{2}$. \vspace{1mm} \label{lem:A-9}
        \item $\displaystyle \frac{\alpha - \beta - |\alpha - \beta|}{2} 
        \leq \alpha - \mathcal{C}(\alpha, \beta) 
        = \alpha - \mathcal{C}(\beta, \alpha) 
        \leq \frac{\alpha - \beta + |\alpha - \beta|}{2}$. \vspace{1mm} \label{lem:A-10}
        \item $|\alpha - \mathcal{C}(\alpha, \beta)| = |\alpha - \mathcal{C}(\beta, \alpha)| \leq |\alpha - \beta|$. \vspace{1mm} \label{lem:A-11}
        \item $\displaystyle \lim_{h_1, h_2 \to 0} \mathcal{C}(\alpha + h_1, -\alpha + h_2) = 0$.  \vspace{1mm}\label{lem:A-12}
        \item $\displaystyle \lim_{\alpha \to \pm \infty}\mathcal{C}(\alpha, \beta) = \beta \pm \sqrt{1 + \beta^2}$ and $\displaystyle \lim_{\beta \to \pm \infty}\mathcal{C}(\alpha, \beta) = \alpha \pm \sqrt{1 + \alpha^2}$.  \label{lem:A-13}
    \end{enumerate}
\end{lemma}

\begin{proof}
    \Cref{lem:A-1,lem:A-2,lem:A-3,lem:A-4} directly follow from the definition of the function $\mathcal{C}$.
    By applying the expression \cref{eq:algebraic_repr_C}, \cref{lem:A-5} can be proved easily.    
    To prove \cref{lem:A-6,lem:A-7}, let $\theta_1 := \arctan\alpha$ and $\theta_2 := \arctan\beta$.
    By \cref{lem:ABC}, $\mathcal{C}(\alpha,\beta) = \tan\left(\frac{\theta_1+\theta_2}{2}\right)$.
    Since the arctangent function is strictly increasing and odd, $\theta_1+\theta_2$ and $\alpha+\beta$ have the same sign.
    Since the tangent function is strictly increasing and vanishes only at zero on $\left(-\frac{\pi}{2},\frac{\pi}{2}\right)$, \cref{lem:A-6,lem:A-7} hold.
    
    Now, we prove \cref{lem:A-8}.
    It is clear that the inequality holds if $\alpha + \beta = 0$.
    Thus, assume that $\alpha + \beta \neq 0$.
    Since
    \begin{displaymath}
        (1+\alpha^2)(1+\beta^2)-(1-\alpha\beta)^2
        =
        (\alpha+\beta)^2
        \geq 0,
    \end{displaymath}
    we have
    \begin{displaymath}
        \sqrt{(1+\alpha^2)(1+\beta^2)}
        \geq
        |1-\alpha\beta|
        \geq
        1-\alpha\beta.
    \end{displaymath}
    Hence,
    \begin{displaymath}
        \alpha\beta-1+\sqrt{(1+\alpha^2)(1+\beta^2)}
        \geq 0.
    \end{displaymath}
    By the Arithmetic Mean-Geometric Mean inequality applied to $\alpha^2 + 1$ and $\beta^2 + 1$, it follows that
    \begin{displaymath}
        0 \leq \alpha \beta - 1 + \sqrt{(1 + \alpha^2)(1 + \beta^2)} \leq \frac{(\alpha + \beta)^2}{2}.
    \end{displaymath}
    If $\alpha + \beta > 0$, dividing by $\alpha + \beta$ yields
    \begin{displaymath}
        \frac{\alpha + \beta - |\alpha + \beta|}{4} =  0 \leq \mathcal{C}(\alpha, \beta) \leq \frac{\alpha + \beta}{2} = \frac{\alpha + \beta + |\alpha + \beta|}{4}.
    \end{displaymath}
    The case $\alpha + \beta < 0$ follows analogously.
    This proves \cref{lem:A-8}.
    Moreover, \cref{lem:A-9} follows immediately from \cref{lem:A-8}.

    Toward proving \cref{lem:A-10}, assume $\alpha + \beta = 0$.
    From the inequality $-|\alpha| \leq 0 \leq |\alpha|$, we obtain
    \begin{displaymath}
        \frac{\alpha + \alpha - |\alpha + \alpha|}{2} \leq \alpha \leq \frac{\alpha + \alpha + |\alpha + \alpha|}{2}.
    \end{displaymath}
    Since $\alpha = - \beta$ and $\mathcal{C}(\alpha, \beta) = 0$, the desired inequalities in \cref{lem:A-10} follow.
    On the other hand, assume that $\alpha + \beta \neq 0$.
    Applying the formula \cref{eq:algebraic_repr_C}, we find that
    \begin{align*}
        \alpha - \mathcal{C}(\alpha, \beta)
        &= \alpha - \frac{\alpha \beta - 1 + \sqrt{(1 + \alpha^2)(1 + \beta^2)}}{\alpha + \beta} 
        = \frac{\sqrt{1 + \alpha^2} \left(\sqrt{1 + \alpha^2} - \sqrt{1 + \beta^2}\right) }{\alpha + \beta}  \\
        &= \frac{\sqrt{1 + \alpha^2} (\alpha^2 - \beta^2) }{(\alpha + \beta)\left(\sqrt{1 + \alpha^2} + \sqrt{1 + \beta^2}\right)} 
        = \frac{\alpha - \beta}{1+ \sqrt{\frac{1 + \beta^2}{1 + \alpha^2}}}.
    \end{align*}
    If $\alpha - \beta \geq 0$, then 
    \begin{displaymath}
        \frac{\alpha - \beta - |\alpha - \beta|}{2} = 0 \leq \frac{\alpha - \beta}{1+ \sqrt{\frac{1 + \beta^2}{1 + \alpha^2}}} \leq \alpha - \beta = \frac{\alpha - \beta + |\alpha - \beta|}{2},
    \end{displaymath}
    and hence the desired inequalities in \cref{lem:A-10} follow. 
    The case when $\alpha - \beta < 0$ can be proved similarly.
    \Cref{lem:A-11} follows from \cref{lem:A-10}.

    Next, to prove \cref{lem:A-12}, apply \cref{lem:A-9} to obtain that
    \begin{displaymath}
        0 \leq \left\vert \, \mathcal{C}(\alpha + h_1, -\alpha + h_2) \right\vert  \leq \frac{|h_1 + h_2|}{2} \leq |h_1| + |h_2| \to 0
    \end{displaymath}
    as $h_1, h_2 \to 0$.

    For \cref{lem:A-13}, observe that
    \begin{displaymath}
        \lim_{\alpha \to \pm \infty} \mathcal{C}(\alpha, \beta)
        =\lim_{\alpha \to \pm \infty} \frac{\beta - \frac{1}{\alpha} \pm \sqrt{\left(1 + \frac{1}{\alpha^2}\right)(1 + \beta^2)}}{1 + \frac{\beta}{\alpha}} 
        = \beta \pm \sqrt{1 + \beta^2}.
    \end{displaymath}
    The case $\beta \to \pm \infty$ can be shown similarly.
\end{proof}

\section{Leicester theorem}
\label{apx:Leicester}

The Leicester theorem, due to W. H. Young, states that, with at most countably many exceptions, the left- and right-hand upper limits coincide, as do the corresponding lower limits \cite[Sect.~6]{2001_Bruckner}.
For completeness, we include a simple direct proof.
Alternatively, the theorem follows immediately from the stronger result known as Young's Rome theorem; for a proof of the latter, see \cite[Thm.~26.2 and Ex.~26.3]{1985_Thomson_BOOK}.

Throughout this section, given an open interval $I\subset\mathbb{R}$, a function $f:I\to\mathbb{R}$, and a point $x\in I$, we write
\begin{displaymath} 
    \underline{f}(x\pm) := \liminf_{h\searrow0}f(x\pm h)
    \qquad\text{and}\qquad
    \overline{f}(x\pm) := \limsup_{h\searrow0}f(x\pm h).
\end{displaymath}

\begin{theorem}[Leicester theorem] \label{thm:Leicester}
    Let $I$ be an open interval in $\mathbb{R}$, and let $f : I \to \mathbb{R}$ be arbitrary.
    Then there exists an at most countable set $N \subset I$ such that, for every $x\in I\setminus N$, $\underline{f}(x+) = \underline{f}(x-)$ and $\overline{f}(x+) = \overline{f}(x-)$.
\end{theorem}

\begin{proof}
    For $p, q\in\mathbb{Q}$ with $p<q$, let 
    \begin{displaymath}
        E_{p,q} := \left\{ x \in I \, \middle| \, \overline{f}(x-) < p < q < \overline{f}(x+) \right\}.
    \end{displaymath}
    
    We claim that for each $x \in E_{p, q}$, there exists $\delta_x > 0$ such that $E_{p, q} \cap (x - \delta_x, x) = \varnothing$.
    Indeed, fix $x \in E_{p, q}$.
    Since $\overline{f}(x-)<p$, there exists $\delta_x > 0$ such that $(x-\delta_x,x)\subset I$ and $f(t) < p$ for every $t\in(x-\delta_x,x)$.
    Suppose, to the contrary, that $E_{p,q}\cap(x-\delta_x,x)\neq\varnothing$, and choose $y \in E_{p,q}\cap(x-\delta_x,x)$.
    Since $q<\overline{f}(y+)$ and $y<x$, the definition of $\overline{f}(y+)$ yields some $z \in (y, x)$ such that $f(z)>q$.
    On the other hand, $(y,x)\subset(x-\delta_x,x)$, and hence $f(z)<p<q$, a contradiction.
    Therefore, $E_{p, q} \cap (x - \delta_x, x) = \varnothing$.
    Since $x \in E_{p, q}$ was arbitrary, the claim holds.

    For $m, n \in\mathbb{N}$, let
    \begin{displaymath}
        F_{m,n}
        :=
        [-m, m] \cap \left\{ x \in E_{p,q} \,\middle|\, E_{p,q}\cap \left(x-\frac{1}{n}, x \right)=\varnothing \right\} .
    \end{displaymath}
    We claim that $F_{m,n}$ is finite.
    Indeed, if $x_1 < x_2 < \cdots < x_k$ are distinct points of $F_{m, n}$, then $x_{j+1} - x_j \geq \frac{1}{n}$ for every $j = 1, 2, \ldots, k - 1$.
    Summing these inequalities over $j=1,2,\ldots,k-1$ gives $\frac{k-1}{n} \leq x_k - x_1 \leq 2m$, and hence $k \leq 2mn+1$.
    Thus the claim follows.

    Moreover,
    \begin{displaymath}
        E_{p,q}
        =
        \bigcup_{m,n\in\mathbb{N}}F_{m,n}.
    \end{displaymath}
    Indeed, for each $x\in E_{p,q}$, one can choose $m,n\in\mathbb{N}$ such that $x \in [-m, m]$ and $\frac{1}{n} < \delta_x$.
    Hence $E_{p,q}$ is at most countable.

    Since
    \begin{displaymath}
        \left\{x\in I\,\middle|\,\overline{f}(x-)<\overline{f}(x+)\right\}
        =
        \bigcup_{\substack{p,q\in\mathbb{Q}\\p<q}}E_{p,q},
    \end{displaymath}
    this set is at most countable.
    Interchanging left and right shows that the set 
    \begin{displaymath}
        \left\{ x \in I \, \middle| \, \overline{f}(x+)<\overline{f}(x-) \right\}
    \end{displaymath}    
    is also at most countable.
    Consequently, $\overline{f}(x-)=\overline{f}(x+)$ outside an at most countable subset of $I$.

    Applying the same argument to $-f$ shows that $\underline{f}(x-)=\underline{f}(x+)$ outside an at most countable subset of $I$.
    Taking the union of the two exceptional sets proves the result.
\end{proof}

We now apply Young's Leicester theorem to establish an auxiliary result needed in the proof of \cref{thm:countable_nonzero_sd_discontinuities}.

\begin{lemma}
    \label{lem:countable_discontinuous_strict_increase}
    Let $I$ be an open interval in $\mathbb{R}$, and let $f:I\to\mathbb{R}$.
    Then $\Disc(f) \cap \bigl(\Inc(f)\cup\Dec(f)\bigr)$ is at most countable.
\end{lemma}

\begin{proof}
    By \cref{thm:Leicester}, there exists an at most countable set $N\subset I$ such that $\underline{f}(x-)=\underline{f}(x+)$ and $\overline{f}(x-)=\overline{f}(x+)$ for every $x \in I\setminus N$.

    Fix $x\in\Inc(f)\setminus N$.
    The definition of $\Inc(f)$ gives $\overline{f}(x-)\leq f(x)\leq \underline{f}(x+)$.
    Hence
    \begin{displaymath}
        \overline{f}(x-)
        \leq f(x)
        \leq \underline{f}(x+)
        \leq \overline{f}(x+)
        =
        \overline{f}(x-).
    \end{displaymath}
    Together with $\underline{f}(x-)=\underline{f}(x+)$, this shows that all four
    one-sided upper and lower limits equal $f(x)$.
    Thus $f$ is continuous at $x$.

    Similarly, every $x\in\Dec(f)\setminus N$ satisfies
    \begin{displaymath}
        \overline{f}(x+)
        \leq f(x)
        \leq \underline{f}(x-)
        =
        \underline{f}(x+)
        \leq \overline{f}(x+),
    \end{displaymath}
    and hence $f$ is again continuous at $x$.

    Combining these inclusions yields $\bigl( \bigl(\Inc(f)\cup\Dec(f)\bigr) \setminus N \bigr) \subset \Cont (f)$.
    Consequently, 
    \begin{displaymath}
        \Disc(f) \cap \bigl(\Inc(f)\cup\Dec(f)\bigr) \subset N.
    \end{displaymath}
    Since $N$ is at most countable, the result follows.
\end{proof}

\bibliographystyle{siamplain}
\bibliography{reference}{}

\end{document}

%% file: shared_info.tex
\usepackage{lipsum}
\usepackage{amsfonts}
\usepackage{amssymb}
\usepackage{bookmark}
\usepackage{bm}
\usepackage{graphicx}
\usepackage{epstopdf}
\usepackage{hyperref}
\usepackage{enumitem}
\usepackage{algorithmic}
\usepackage{booktabs}
\usepackage{multirow}
\usepackage{cleveref}
\usepackage{aliascnt}
\usepackage[caption=false]{subfig}
\usepackage{orcidlink}

\ifpdf
  \DeclareGraphicsExtensions{.eps,.pdf,.png,.jpg}
\else
  \DeclareGraphicsExtensions{.eps}
\fi

\newsiamremark{remark}{Remark}
\newsiamremark{hypothesis}{Hypothesis}
\crefname{hypothesis}{Hypothesis}{Hypotheses}
\newsiamthm{claim}{Claim}
\newsiamremark{fact}{Fact}
\crefname{fact}{Fact}{Facts}

\newaliascnt{example}{theorem}
\newtheorem{example}[example]{\textit{Example}}
\aliascntresetthe{example}
\crefname{example}{Example}{Examples}
\Crefname{example}{Example}{Examples}
\crefname{remark}{Remark}{Remarks}
\Crefname{remark}{Remark}{Remarks}

\crefname{enumi}{part}{parts}
\Crefname{enumi}{Part}{Parts}

\newcommand{\sd}{\mathord{\prime\mkern-2.5mu\reflectbox{$\scriptstyle\prime$}}}

\Crefname{ALC@unique}{Line}{Lines}

\headers{Specular differentiation}{K. Jung}

\title{Specular differentiation in one dimension: a quasi-mean value theorem, regularity, and discontinuities\thanks{Date: \today.
\funding{This work is supported in part by Michigan State University and the National Science Foundation Research Traineeship Program (DGE-2152014). The author received partial support from Jun Kitagawa's National Science Foundation grant (DMS-2246606).}}
}

\author{Kiyuob Jung\,\orcidlink{0009-0007-5075-4061}\thanks{Department of Mathematics, Michigan State University, East Lansing, MI 48824 USA (\email{kyjung@msu.edu}).}
}
\usepackage{amsopn}
\DeclareMathOperator{\sgn}{sgn}
\DeclareMathOperator{\Disc}{Disc}
\DeclareMathOperator{\Cont}{Cont}
\DeclareMathOperator{\Inc}{Inc}
\DeclareMathOperator{\Dec}{Dec}
\DeclareMathOperator{\Span}{span}